\documentclass{article}

\usepackage{amssymb}
\usepackage{amsmath,amsfonts,amsthm}
\usepackage[hmargin=2.5cm,vmargin=3.5cm,bmargin=2.5cm]{geometry}
\usepackage{graphicx}
\graphicspath{ {images/} }

\newtheorem{theorem}{Theorem}[section]
\newtheorem{lemma}[theorem]{Lemma}
\newtheorem{proposition}[theorem]{Proposition}

\newcommand{\cqd}{\hfill $\Box$ \vspace{0.2cm}}

\begin{document}


\title{Third order differential equations and local isometric immersions of pseudospherical surfaces}

\author{Tarc\'isio Castro Silva \qquad Niky Kamran}

\date{}

\maketitle{}

\begin{abstract}
The class of differential equations describing pseudospherical surfaces enjoys important integrability properties which manifest themselves by the existence of infinite hierarchies of conservation laws (both local and non-local) and the presence associated linear problems. It thus contains many important known examples of integrable equations, like the sine-Gordon, Liouville, KdV, mKdV, Camassa-Holm and Degasperis-Procesi equations, and is also home to many new families of integrable equations. Our paper is concerned with the question of the local isometric immersion  in ${\bf E}^{3}$ of the pseudospherical surfaces defined by the solutions of equations belonging to the class of Chern and Tenenblat~\cite{CT}. In the case of the sine-Gordon equation, it is a classical result that the second fundamental form of the immersion depends only on a jet of finite order of the solution of the pde. A natural question is therefore to know if this remarkable property extends to equations other than the sine-Gordon equation within the class of differential equations describing pseudospherical surfaces. In a pair of earlier papers~\cite{KKT1},~\cite{KKT2} we have shown that this property fails to hold for all $k$-th order evolution equations $u_t= F(u,u_x,..., u_{x^k})$ and all other second order equations of the form $u_{xt}=F(u,u_{x})$, except for the sine-Gordon equation and a special class of equations for which the coefficients of the second fundamental form are universal, that is functions of $x$ and $t$ which are independent of the choice of solution $u$. In the present paper, we consider third-order equations of the form $u_{t}-u_{xxt}=\lambda u u_{xxx} + G(u,u_x,u_{xx}),\, \lambda \in \mathbb{R},
$ which describe pseudospherical surfaces. This class contains the Camassa-Holm and Degasperis-Procesi equations as special cases. We show that whenever there exists a local isometric immersion in ${\bf E}^3$ for which the coefficients of the second fundamental form  depend on a jet of finite order of $u$, then these coefficients are universal in the sense of being independent on the choice of solution $u$. This result further underscores the special place that the sine-Gordon equations seems to occupy amongst integrable partial differential equations in one space variable.

\vspace{0.3 cm}

\end{abstract}

\section{Introduction}
The class of partial differential equations describing pseudospherical surfaces, first introduced in a fundamental paper by Chern and Tenenblat~\cite{CT}, gives a rich geometric framework for the classification and study of integrable partial differential equations in one space variable. Recall that a partial differential equation 
\begin{equation}\label{pde}
\Delta(t,x,u,u_{t},u_{x},\ldots,u_{t^{l}x^{k-l}})=0,
\end{equation}
is said to describe pseudospherical surfaces if there exist $1$-forms 
\begin{equation}\label{forms}
\omega_{i}=f_{i1}dx+f_{i2}dt,\quad 1\leq i \leq3,
\end{equation}
where the coefficients $f_{ij},\,1\leq i \leq 3,\,1\leq j\leq 2,$ are smooth functions of $t,x,u$ and finitely many derivatives of $u$, such that the structure equations of a surface of Gaussian curvature equal to $-1$,
\begin{equation}\label{struct}
d\omega_{1}=\omega_{3} \wedge\omega_{2},\quad d\omega_{2}=\omega_{1} \wedge\omega_{3},\quad d\omega_{3}=\omega_{1} \wedge\omega_{2},
\end{equation}
are satisfied if and only if $u$ is a smooth solution of (\ref{pde}) for which 
\begin{equation}\label{indep}
\omega_{1}\wedge \omega_{2}\neq 0.
\end{equation}
It thus follows that every such solution determines a pseudospherical metric, that is a Riemannian metric of constant negative Gaussian curvature equal to $-1$, defined by
\begin{equation}\label{metric} 
ds^{2}=(\omega_{1})^{2}+(\omega_{2})^{2}.
\end{equation}
The 1-form $\omega_{3}$ appearing in the structure equations~(\ref{struct}) is then the Levi-Civita connection $1$-form of the metric~(\ref{metric}).

The prototypical example of a differential equation describing pseudospherical surfaces is the sine-Gordon equation
\begin{equation}\label{sG}
u_{tx}=\sin u, 
\end{equation}
for which a choice of $1$-forms (\ref{forms}) satisfying the structure equations (\ref{struct}) is given by
\begin{equation}\label{sGcof1}
\omega_1 = \frac{1}{\eta}\sin u \, dt, \quad \omega_2 = \eta\, dx+\frac{1}{\eta}\cos u \,dt,\quad \omega_3 = u_{x}\,dx.
\end{equation}
There may of course be different choices of $1$-forms satisfying the structure equations (\ref{struct}) for a given differential equation describing pseudo spherical surfaces; for example, for the sine-Gordon equation (\ref{sG}), a choice different from the one given in (\ref{sGcof1}) would be
\begin{equation}\label{sGcof2}
\omega_1 = \cos \frac{u}{2}( dx+dt),\quad \omega_2 =\sin \frac{u}{2} (dx - dt),\quad \omega_3 = \frac{u_{x}}{2} dx - \frac{u_{t}}{2} dt.
 \end{equation}
In~(\ref{sGcof1}), the constant $\eta$ is a continuous non-zero real parameter which reflects the existence of a one-parameter family of B\"acklund transformation for the sine-Gordon equation and is key to the existence of infinitely many conservation laws. More generally one may consider partial differential equations describing pseudospherical surfaces with the property that one of the components $f_{ij}$ can be chosen to be a continuous parameter are said to describe $\eta$ pseudospherical surfaces. Each equation belonging to this class is the integrability condition of a linear system of the form 
\begin{equation}\label{linear}
dv^{1}=\frac{1}{2}(\omega_2{}v^{1}+(\omega_{1}-\omega_{3})v^{2}),\quad dv^{2}=\frac{1}{2} ((\omega_{1}+\omega_{3})v^{1}-\omega_{2}v^{2}),
\end{equation}
which may be used to solve the differential equation by inverse scattering~\cite{BRT}, with $\eta$ playing the role of a spectral parameter for the scattering problem. It is also shown in~\cite{CaT} that one can generate infinite sequences of conservation laws for the class of differential equations describing $\eta$ pseudospherical surfaces by making use of the structure equations~(\ref{struct}), although some of these conservation laws may end up being non-local. Important further developments of these ideas around this theme can be found in~\cite{R88},~ \cite{R89},~ \cite{RT90}, ~\cite{RT92},~ \cite{KT}, ~ \cite{R98}, ~\cite{reyes},~\cite{FT},~\cite{Go},~\cite{GR},~\cite{JT},~\cite{FOR},~ \cite{R106}.

One may also consider the class of differential equations describing pseudospherical surfaces from an extrinsic point of view, motivated by the classical result which says that every pseudospherical surface can be locally isometrically immersed in $\bf{E}^{3}$. One would expect the dependence of the second fundamental form of the immersion on the solution chosen for the differential equation to be quite complicated. However, the formula for the second fundamental form turns out to be particularly simple in the case of the sine-Gordon equation, as we now recall. Indeed, we first recall that the components $a,b,c$ of the second fundamental form of a local isometric immersion of a pseudospherical surface into $\bf{E}^{3}$ are defined by the relations
\begin{equation}\label{3.4}
\omega_{13} = a\omega_1+b\omega_2, \quad \omega_{23} = b\omega_1+c\omega_2,
\end{equation}
where the $1$-forms $\omega_{13}, \omega_{23}$ satisfy the structure equations
\begin{equation}\label{3.5}
d\omega_{13} = \omega_{12}\wedge\omega_{23}, \quad d\omega_{23} = \omega_{21}\wedge\omega_{13},
\end{equation}
equivalent to the Codazzi equations, and the Gauss equation for a pseudo spherical surface, given by
\begin{equation}\label{gauss}
ac-b^2=-1.
\end{equation}

We recall from \cite{KKT1} that the Codazzi equations \eqref{3.5} may be expressed in terms of the components $f_{ij}$ of the 1-forms $\omega_1$, $\omega_2$, $\omega_3$ in the following form
\begin{eqnarray}
f_{11}D_t a+f_{21}D_t b -f_{12}D_x a-f_{22}D_x b-2b\Delta_{13}+(a-c)\Delta_{23} = 0,&&\label{3.8}\\
f_{11}D_t b+f_{21}D_t c -f_{12}D_x b-f_{22}D_x c+(a-c)\Delta_{13}+2b\Delta_{23} = 0,&&\label{3.9}
\end{eqnarray}
where
\begin{eqnarray}\label{3.7}
\Delta_{ij} = f_{i1}f_{j2} - f_{j1}f_{i2}.
\end{eqnarray}
and where we assume~\cite{KKT1} that 
\begin{eqnarray}\label{3.10a}
\quad \Delta_{13}^2+\Delta_{23}^2 \neq 0.
\end{eqnarray}

For the sine-Gordon equation, with the choice of $1$-forms $\omega_{1}, \omega_{2}$ and $\omega_{3}=\omega_{12}$ given by (\ref{sGcof2}), it is easily verified that the $1$-forms $\omega_{13},\omega_{23}$ are given by
\begin{eqnarray*}
\omega_{13} &=& \sin\frac{u}{2} (dx+dt) =  \tan \frac{u}{2}\omega_1,\\
\omega_{23} &=& -\cos\frac{u}{2} (dx - dt) = -\cot \frac{u}{2}\omega_2.
\end{eqnarray*}
What is particularly noteworthy in the case of the sine-Gordon equation that the components $a,b,c$ depend only on $u$ and {\emph{finitely many derivatives} of $u$. It is therefore a natural question to ask whether such a remarkable property holds for other equations within the class of differential equations describing pseudospherical surfaces, or whether the sine-Gordon equation in any way special in this regard. In~\cite{KKT1} and~\cite{KKT2}, we investigated this question for  $k$-th order evolution equations
\begin{equation}\label{evoleq}
u_t= F(u,u_x,..., u_{x^k}),
\end{equation}
and second order hyperbolic equations
\begin{equation}\label{fGord}
u_{xt}=F(u,u_{x}),
\end{equation}
and proved that there are no other equations than the sine-Gordon equation for which this property holds, except for some special equations for which $a,b,c$ are {\emph{universal}}, that is functions of $x,t$ which are {\emph{independent of}} $u$. These results show that the sine-Gordon equation occupies a special position within the class of differential equations of the form~(\ref{evoleq}) and~(\ref{fGord}) which describe pseudospherical surfaces.

Our goal in the present paper is to investigate this question for the class of partial differential equations given by 
\begin{eqnarray}\label{T}
u_{t}-u_{xxt}=\lambda u u_{xxx} + G(u,u_x,u_{xx}),\quad \lambda \in \mathbb{R},
\end{eqnarray}
which describe pseudospherical surfaces under the condition of that 1-forms $\omega_i=f_{i1}dx+f_{i2}dt$ satisfy
\begin{equation}\label{1}
f_{p1}=\mu_{p}f_{11}+\eta_{p}, \quad \mu_p, \eta_p \in \mathbb{R}, \quad 2\leq p \leq 3.
\end{equation}
This class of equations, which has recently been classified by Castro Silva and Tenenblat~\cite{CST} contains important examples such as the Camassa-Holm equation \cite{CH}
\begin{eqnarray*}
u_{t}-u_{xxt} = u u_{xxx} + 2 u_x u_{xx} - 3 u u_x-m u_x, \quad m \in \mathbb{R},
\end{eqnarray*}
and Degasperis-Procesi equation \cite{DP}
\begin{eqnarray*}
u_{t}-u_{xxt} = u u_{xxx} + 3 u_x u_{xx} - 4 u u_x.
\end{eqnarray*}
Our main result is the following:

\begin{theorem}\label{teo}
Except for the two families of third order differential equations of the form
\begin{eqnarray}\label{1a}
u_t - u_{xxt} = \frac{1}{h'}(m\psi + \psi_x), \quad m\in \mathbb{R}\setminus \left\lbrace 0\right\rbrace,
\end{eqnarray}
where $\psi(u,u_x)\neq 0$ and $h(u-u_{xx})$ are differentiable functions, with $h'\neq 0$, and
\begin{eqnarray}\label{2a}
u_t - u_{xxt} = \lambda uu_{xxx}+\frac{1}{h'}[m_1\psi + \psi_x -\lambda uu_xh' - (\lambda u_x+\lambda m_1u+m_2)h], 
\end{eqnarray}
where $\lambda$, $m_1$, $m_2$ $\in$ $\mathbb{R}$, $(\lambda m_1)^2+m_2^2\neq 0$, $\psi(u,u_x)$ and $h(u-u_{xx})$ are differentiable functions, with $h'\neq 0$, there exists no third order partial differential equation of type \eqref{T} describing pseudospherical surfaces, under the condition \eqref{1}, with the property that the coefficients of the second fundamental forms of the local isometric immersions of the surfaces associated to the solutions $u$ of the equation depend on a jet of finite order of u. Moreover, the coefficients of the second fundamental forms of the local isometric immersions of the surfaces determined by the solution $u$ of \eqref{1a} or \eqref{2a} are universal, i.e., they are universal functions of $x$ and $t$, independent of $u$.

\end{theorem}

We see in particular the Degasperis-Procesi equation belongs to the class \eqref{2a} of equations covered by Theorem~\ref{teo}. On the other hand, the Camassa-Holm equation is not covered by either \eqref{1a} or \eqref{2a}, meaning that for the Degasperis-Procesi equation, the components $a,b,c$ of the second fundamental form are the same universal functions of $x$ and $t$ for any solution $u$, while for the Camass-Holm equation the components $a,b,c$ depend on jets of arbitrary high order of $u$. Theorem~\ref{teo} underscores once again the special place that the sine-Gordon equations appears to occupy amongst integrable partial differential equations in one space variable.

Our paper is organized as follows. In Section~\ref{sec2}, we recall without proof the classification results of~\cite{CST} that will be needed to prove Theorem~\ref{teo}. The classification splits into branches which are treated on a case-by-case basis in Section~\ref{sec4}, starting from the expression of the Codazzi and Gauss equations in terms of the coefficients $f_{ij}$ of the $1$-forms $\omega_{1},\,\omega_{2},\,\omega_{3}$ (see \eqref{3.8}, \eqref{3.9}). Finally, we carry out in Section~\ref{sec5} the integration of the Codazzi and Gauss equations in the cases in which the components $a,b,c$ of the second fundamental form are universal functions of $x$ and $t$ and obtain explicit expressions for these functions.

\section{Third order differential equations describing pseudospherical surfaces}\label{sec2}

\hspace{0.5 cm}Let us recall from \cite{CST} without proof the characterization and classification theorems (Theorems \ref{teo7.1}, \ref{teo7.2}-\ref{teo7.5}) of the equations \eqref{T} that describe pseudospherical surfaces under the hypothesis \eqref{1}. We will use the following notation, also used in~\cite{CT}, for the spatial derivatives of $u$,
$$
z_i = \partial_x^i u, \quad 0\leq i
$$

\begin{theorem}\label{teo7.1}
\textnormal{\cite{CST}} An equation
\begin{eqnarray}\label{eqT}
z_{0,t}-z_{2,t}=\lambda z_0z_3 + G(z_0,z_1,z_2), \quad G\neq 0,
\end{eqnarray}
describes pseudospherical surfaces, with associated 1-forms $\omega_i=f_{i1}dx+f_{i2}dt$, $1 \leq i \leq 3$, where $f_{ij}$ are real and differentiable functions of $z_k$, $0\leq k \leq \ell$, $\ell$ $\in$ $\mathbb{Z}$, satisfying \eqref{1} if, and only if, $f_{ij}$ and $G$ satisfy
\begin{eqnarray}\label{2.3.1}
f_{11,z_0}\neq 0, \quad f_{11,z_0} + f_{11,z_2} = 0, \quad f_{i2,z_s} = 0, \quad f_{11,z_1} = f_{11,z_s} = 0,\quad s \geq 3,
\end{eqnarray}
\begin{eqnarray}\label{(7.2)}
f_{i2} = -\lambda z_0 f_{i1} + \phi_{i2},
\end{eqnarray}
where $\phi_{i2}(z_0,z_1)$ are real and differentiable functions of $z_0$ and $z_1$ satisfying
\begin{eqnarray}
-f_{11,z_0}G+\displaystyle \sum_{i=0}^{1} z_{i+1}f_{12,z_i} + (\mu_2 \phi_{32}-\mu_3 \phi_{22})f_{11} + \eta_2 \phi_{32}-\eta_3 \phi_{22}=0,&&\label{(7.3)}\\
-\mu_2 f_{11,z_0}G+\displaystyle \sum_{i=0}^{1} z_{i+1}f_{22,z_i} - (\phi_{32}-\mu_3 \phi_{12})f_{11} + \eta_3 \phi_{12}=0,&&\label{(7.4)}\\
-\mu_3 f_{11,z_0}G+\displaystyle \sum_{i=0}^{1} z_{i+1}f_{32,z_i} - (\phi_{22}-\mu_2 \phi_{12})f_{11} + \delta \eta_2 \phi_{12}=0,&&\label{(7.5)}\\
(\phi_{22}-\mu_2 \phi_{12})f_{11} - \eta_2 \phi_{12}\neq 0.&&\label{7.5.1}
\end{eqnarray}
\end{theorem}


\begin{theorem}\label{teo7.2}\textnormal{\cite{CST}} Consider an equation of type $(\ref{T})$ that describes pseudospherical surfaces, with associated 1-forms $\omega_i=f_{i1}dx+f_{i2}dt$, $1 \leq i \leq 3$, where $f_{ij}$ are real and differentiable functions of $z_k$, satisfying \eqref{1} and \eqref{2.3.1}-\eqref{7.5.1}. Suppose that $\phi_{22} - \mu_2 \phi_{12}\equiv 0$ and $\mu_2\mu_3 \eta_2 -(1+ \mu_2^2)\eta_3=0$. Then the equation is given by 
\begin{eqnarray}\label{th2.2}
z_{0,t}-z_{2,t} = \frac{1}{h'}(z_1\psi_{,z_0} + z_2 \psi_{,z_1} + m\psi), \quad m \in \mathbb{R}, \quad m\neq 0,
\end{eqnarray}
\begin{eqnarray}\label{th2.2a}
 \begin{array}{ll}
& f_{11}=h,\\ & f_{21}=\mu h \pm m\sqrt{1+\mu^2},\\ & f_{31}=\pm \sqrt{1+\mu^2} h + m\mu, 
\end{array}\quad
\begin{array}{ll}
f_{12}=\psi, \\ f_{22}=\mu \psi, \\ f_{32} = \pm \sqrt{1+\mu^2} \psi,
\end{array}
\end{eqnarray}
where $\lambda = 0$ and $\mu \in \mathbb{R}$,  $h(z_0-z_2)$ and $\psi(z_0,z_1)$ are real and differentiable functions satisfying $h'\neq 0$ and $\psi \neq 0$.
\end{theorem}

\begin{theorem}\label{teo7.3}\textnormal{\cite{CST}} 
Consider an equation of type $(\ref{T})$ that describes pseudospherical surfaces, with associated 1-forms $\omega_i=f_{i1}dx+f_{i2}dt$, $1 \leq i \leq 3$, where $f_{ij}$ are real and differentiable functions of $z_k$, satisfying \eqref{1} and \eqref{2.3.1}-\eqref{7.5.1}. 
Suppose that $\phi_{22} - \mu_2 \phi_{12}\equiv 0$ and $\mu_2\mu_3 \eta_2 -(1+ \mu_2^2)\eta_3\neq 0$. Then the equation is given by 
\begin{eqnarray}\label{th2.3}
z_{0,t}-z_{2,t} =\lambda z_0z_3 -\frac{\lambda}{h'}(z_1 h + z_0z_1 h' + m_1z_1 + m_2z_2), \quad \lambda, \textnormal{ } m_1,\textnormal{ } m_2 \in \mathbb{R}, \quad \lambda m_2\neq 0,
\end{eqnarray}
\begin{eqnarray}\label{th2.3a}
\begin{array}{ll}
 f_{11}=h,\\  f_{21} = \mu h + \eta,\\  f_{31} = \left[\frac{m_1(1+\mu^2)}{m_2\eta}-\frac{\mu}{m_2}\right]h+\frac{m_1\mu - \eta}{m_2}, 
\end{array}
\begin{array}{ll}
f_{12} = -\lambda z_0 h - \lambda m_2 z_1, \\ f_{22} = -\lambda \mu z_0 h - \lambda m_2\mu z_1 - \lambda \eta z_0, \\ f_{32} = -\lambda z_0f_{31} - \frac{\lambda}{\eta}\left[m_1(1+\mu^2)-\mu\eta\right]z_1,
\end{array}
\end{eqnarray}
where $\mu,\eta\in \mathbb{R}$, $\eta \neq 0$ and  $h(z_0-z_2)$ is a real differential function of $z_0-z_2$ satisfying $h'\neq 0$ and
\[
(m_2\eta)^2=m_1^2 + (m_1\mu-\eta)^2.
\]
\end{theorem}

\begin{theorem}\label{teo7.4} \textnormal{\cite{CST}}
Consider an equation of type $(\ref{T})$ that describes pseudospherical surfaces, with associated 1-forms $\omega_i=f_{i1}dx+f_{i2}dt$, $1 \leq i \leq 3$, where $f_{ij}$ are real and differentiable functions of $z_k$, satisfying \eqref{1} and \eqref{2.3.1}-\eqref{7.5.1}. 
Suppose that $\phi_{22} - \mu_2 \phi_{12}\neq 0$ and $\mu_2\mu_3 \eta_2 -(1+ \mu_2^2)\eta_3 = 0$. Then the equation is given by 
\begin{eqnarray}\label{th2.4}
z_{0,t}-z_{2,t} =\lambda z_0z_3+ \frac{1}{h'}\left[z_2\psi_{,z_1}+z_1\psi_{,z_0} +  m_1\psi
- \lambda z_0z_1h'
-\left(\lambda z_1 + \lambda m_1z_0 + m_2\right) h  \right],
\end{eqnarray}
\begin{eqnarray}\label{th2.4a}
\begin{array}{ll}
& f_{11}=h,\\ & f_{21} = \mu h \pm m_1\sqrt{1+\mu^2},\\ & f_{31} = \pm \sqrt{1+\mu^2} h +m_1\mu, 
\end{array}
\begin{array}{ll}
f_{12} = -\lambda z_0 h + \psi, \\ f_{22} = -\lambda \mu z_0 h + \mu \psi \pm m_2\sqrt{1+\mu^2}, \\ 
f_{32} = \pm \sqrt{1+\mu^2}(\psi- \lambda z_0 h)  +\mu m_2,
\end{array}
\end{eqnarray}
where $\mu,\, \lambda,\,m_1,\, m_2$ $\in \mathbb{R}$, $(\lambda m_1)^2+m_2^2 \neq 0$,  $\,h(z_0 - z_2)$ and $\psi(z_0,z_1)$ are real and differentiable functions, with $h'\neq 0$.
\end{theorem}

\begin{theorem}\label{teo7.5}\textnormal{\cite{CST}}  
Consider an equation of type $(\ref{T})$ that describes pseudospherical surfaces, with associated 1-forms $\omega_i=f_{i1}dx+f_{i2}dt$, $1 \leq i \leq 3$, where $f_{ij}$ are real and differentiable functions of $z_k$, satisfying \eqref{1} and \eqref{2.3.1}-\eqref{7.5.1}. 
Suppose that $\phi_{22} - \mu_2 \phi_{12}\neq 0$ and $\mu_2\mu_3 \eta_2 -(1+ \mu_2^2)\eta_3 \neq 0$. Then the equation is given by 

\begin{flushleft}
$\left( i \right)$ 
\hspace{0.3 cm}$ z_{0,t}-z_{2,t}= \lambda z_0z_3+ \lambda \left(z_1z_2 - 2 z_0z_1 - \frac{m}{\tau}z_1 \mp  \frac{z_2}{\tau}\right)
+ \tau e^{\pm \tau z_1}\left( \tau z_0z_2 \pm z_1+ m z_2\right)\varphi$
\begin{eqnarray}\label{th2.5i}
\hspace*{2cm}\pm e^{\pm \tau z_1}\left(\tau z_0z_1 +\tau z_1z_2 + mz_1\pm z_2\right)\varphi'
+ z_1^2 e^{\pm \tau z_1}\varphi'',\quad  \lambda, m, \tau \in \mathbb{R}, \quad \tau>0,
\end{eqnarray}
\begin{eqnarray}\label{th2.5a}
\begin{array}{ll}
f_{11}=a(z_0-z_2)+b,\\
f_{21} =  \mu f_{11}+\eta,\\
f_{31} = \pm \left(m-\frac{b\tau}{a}\right)\left(\frac{1+\mu^2}{\eta}f_{11}+\mu\right)\mp \frac{\tau}{a}f_{21},\\ 
f_{12} =-\lambda z_0f_{11} + \left[\pm\tau(az_0+b)\varphi + az_1\varphi'\right]e^{\pm\tau z_1}\mp\frac{\lambda a}{\tau}z_1,\\
f_{22} = \mu f_{12}- \lambda\eta z_0 \pm \eta\tau e^{\pm\tau z_1}\varphi,\\ 
f_{32} = \pm\left(m-\frac{b\tau}{a}\right)\left[\frac{1+\mu^2}{\eta}f_{12}-\mu(\lambda z_0\mp\tau e^{\pm\tau z_1}\varphi)\right]\mp\frac{\tau}{a}f_{22},
\end{array}
\end{eqnarray}

\noindent
where $\mu,\, \eta,\, a,\, b$ $\in$ $\mathbb{R}$, $a\eta\neq 0$, $\varphi(z_0)\neq 0$ is a real differentiable function and
\[
(a\eta)^2=(am-b\tau)^2+[\mu(am-b\tau)-\tau\eta]^2.
\]
\vspace*{0.5 cm}
or
\vspace*{0.3 cm}

$\left( ii \right)$ 
\hspace{0.57 cm}
\begin{eqnarray}\label{th2.5ii}
z_{0,t}-z_{2,t}=\lambda z_0z_3+\lambda(2z_1z_2-3z_0z_1-m_2z_1) +
m_1\theta e^{\theta z_0}(\theta z_1^3+z_1z_2+2z_0z_1+m_2z_1)
\end{eqnarray}


\noindent with $\lambda$, $\theta$, $m_1, \, m_2 \in \mathbb{R}$,  $\theta\neq 0$, $\lambda^2+m_1^2\neq 0,$ 
\end{flushleft}
\begin{eqnarray}\label{th2.5b}
\begin{array}{l}
f_{11} = a(z_0-z_2)+b,\\
f_{21} = \mu f_{11}+\eta,\\
f_{31} = \pm \sqrt{1+\mu^2}f_{11} \pm  \frac{\theta+a\mu\eta}{a\sqrt{1+\mu^2}},\\
f_{12} = -\lambda z_0f_{11}+ am_1\theta e^{\theta  z_0}z_1^2+(m_1\theta e^{\theta z_0}-\lambda) 
\left[\frac{az_0+b}{\theta} \pm (\mu- \frac{a\eta}{\theta})\frac{z_1}{\sqrt{1+\mu^2}}  \right],\\

f_{22} = -\lambda z_0f_{21}+\mu a  m_1\theta e^{\theta z_0}z_1^2
+ (m_1\theta e^{\theta z_0}-\lambda) \frac{1}{\theta}
\left[\mu(az_0+b)+\eta  \mp (\theta +\mu \eta a)\frac{z_1}{\sqrt{1+\mu^2}}
  \right],\\

f_{32} = -\lambda z_0 f_{31} \pm \sqrt{1+\mu^2} a m_1\theta e^{\theta z_0}z_1^2
-(m_1\theta e^{\theta z_0}-\lambda) \frac{1}{\theta}
\left\{a\eta z_1 \mp\frac{1}{\sqrt{1+\mu^2}}
\left[(1+\mu^2)(az_0+b)+\mu\eta +\frac{\theta}{a}\right]\right\} 
\end{array}
\end{eqnarray}

\noindent
where $\mu,\, \eta,\, a \in \mathbb{R}$, $a\neq 0$ and  
\begin{equation*}
b=\frac{a}{2\theta}\left[ 
\frac{(\mu\theta-\eta a)^2}{a^2(1+\mu^2)}-\frac{a}{\theta}+m_2\theta-1\right].
\label{relaigual}
\end{equation*}
\end{theorem}



\section{Proof of Theorem \ref{teo}}\label{sec4}

\subsection{Total derivatives and prolongations}

Let us first introduce a compact notation for the time derivatives and mixed derivatives of $u$ in addition to the notation introduced earlier for the spatial derivatives of $u$, by letting
\begin{eqnarray}\label{3.10}
z_i = \partial_x^i u, \quad w_j = \partial_t^j u, \quad v_k = \partial_t^k u_x,
\end{eqnarray}
where $z_0 = w_0 = u$ and $z_1 = v_0 = u_x$. We have therefore,
\begin{eqnarray*}
\begin{array}{ll}
z_{i,x} = z_{i+1}, \\ z_{i,t} = \partial_x^{i-2} u_{xxt},
\end{array}\quad
\begin{array}{ll}
w_{j,x} =  \partial_t^j u_x, \\ w_{j,t} = w_{j+1},
\end{array}\quad
\begin{array}{ll}
v_{k,x}=\partial_t^{k-1}u_{xxt}, \\ v_{k,t}=v_{k+1},
\end{array}\quad
\end{eqnarray*}
and the total derivatives of a differentiable function $\phi = \phi(x,t,z_0,z_1, w_1, v_1, ..., z_{l}, w_{m}, v_{n})$, where $1\leq l  <\infty$, $1\leq m < \infty$ and $1\leq n < \infty$ are finite, but otherwise arbitrary, are given by
\begin{eqnarray}
&& D_x \phi = \phi_x + \sum_{i=0}^l \phi_{z_i}z_{i+1}+ \sum_{j=1}^{m} \phi_{w_j}w_{j,x}+ \sum_{k=1}^{n} \phi_{v_k}v_{k,x},\label{3.11}\\
&& D_t \phi = \phi_t + \sum_{i=2}^l \phi_{z_i}z_{i,t}+ \sum_{j=0}^{m} \phi_{w_j}w_{j+1}+ \sum_{k=0}^{n} \phi_{v_k}v_{k+1}.\label{3.12}
\end{eqnarray}
In particular, we obtain the following expressions for the prolongations of the partial differential equation \eqref{T}
\begin{eqnarray}\label{4.1.1}
z_{2q,t} = z_{0,t} - \sum_{i=0}^{q -1}D_x^{2i}F, \quad z_{2q +1,t} = z_{1,t} - \sum_{i=0}^{q -1}D_x^{2i+1}F,
\end{eqnarray}
where $q=1,2,3,\ldots$, $F(z_0,z_1,z_2,z_3)=\lambda z_0z_3+G(z_0,z_1,z_2)$ and $D_x^0 F = F$.

\subsection{Necessary conditions for the existence of second fundamental forms depending on jets of finite order of $u$}

Our goal in this section is to analyze the system \eqref{gauss}, \eqref{3.8}, \eqref{3.9} governing the components $a$, $b$, $c$ of the second fundamental form and to obtain necessary conditions for the existence of solution depending on jets of finite order of $u$. We note that since the coefficients $f_{ij}$ appearing in the classification given in the Section \ref{sec2}, i.e., in Theorems \ref{teo7.2}-\ref{teo7.5}, depend only on $z_0$, $z_1$ and $z_2$, it follows that the functions $\Delta_{ij}$ defined in \eqref{3.7} depend only on $z_0$, $z_1$ and $z_2$.

\begin{lemma}\label{lemma3.1}
Consider an equation of type \eqref{T} describing pseudospherical surfaces, under the condition \eqref{1}, given by the Theorems $\ref{teo7.2}$-$\ref{teo7.5}$. Assume there is a local isometric immersion of the pseudospherical surface, determined by a solution $u(x,t)$ of \eqref{T} satisfying ~\eqref{indep}, for which the coefficients $a$, $b$ and $c$ of the second fundamental form depend on $x, t, z_0, \ldots, z_l, w_1, \ldots, w_m, v_1, \ldots, v_n$, where $1\leq l  <\infty$, $1\leq m < \infty$ and $1\leq n < \infty$ are finite, but otherwise arbitrary. Then $ac\neq 0$ on any open set of the domain of $u$.
\end{lemma}

\noindent\textbf{Proof}. Firstly, we will show that $c$ is not zero. Then, using the fact $c\neq 0$, we will show that $a=0$ leads to a contradiction and, thus, conclude that $ac\neq 0$.

Assume $c=0$ on a open set. Then, \eqref{gauss} implies $b=\pm 1$ and \eqref{3.8} and \eqref{3.9} reduce to
\begin{eqnarray}
f_{11}D_ta - f_{12}D_xa \mp 2\Delta_{13}+a\Delta_{23} = 0,&&\label{3.14a}\\
a\Delta_{13}\pm 2\Delta_{23} = 0.&&\label{3.15a}
\end{eqnarray}
It follows from \eqref{3.15a} and \eqref{3.10a} that $\Delta_{13}\neq 0$ and $a = \mp 2\Delta_{23}/ \Delta_{13}$. Since $\Delta_{13}$ and $\Delta_{23}$ depend only on $z_0$, $z_1$ and $z_2$, we conclude that $a$ depends only on $z_0$, $z_1$ and $z_2$ and \eqref{3.14a} reduces to
\begin{eqnarray}\label{3.16a}
f_{11}[(a_{z_0}+a_{z_2})z_{0,t}+a_{z_1}z_{1,t}-a_{z_2}(\lambda z_0z_3+G)]-f_{12}\sum_{i=0}^2 a_{z_i}z_{i+1}\mp 2\Delta_{13}+a\Delta_{23} = 0.
\end{eqnarray}
Differentiation with respect to $z_{0,t}$, $z_{1,t}$ and $z_3$ implies
\begin{eqnarray}\label{3.17a}
f_{11}a_{z_1}=f_{11}(a_{z_0}+a_{z_2}) = a_{z_2}(f_{12}+\lambda z_0f_{11}) = 0,
\end{eqnarray}
where we recall from \eqref{(7.2)} that $f_{12}+\lambda z_0f_{11}=\phi_{12}$. Since $f_{11}=h$ can not be zero on any open set (see \eqref{2.3.1}), we have $a_{z_1}=a_{z_0}+a_{z_2} = 0$.

If $\phi_{12}\neq 0$ then from \eqref{3.17a} we conclude that $a$ is a constant and \eqref{3.16a} reduces to $\mp 2\Delta_{13}+a\Delta_{23} = 0$. This equation with \eqref{3.15a} implies that $\Delta_{13}=\Delta_{23} = 0$ which contradicts \eqref{3.10a}.

If $\phi_{12}=0$ on a open set, the only equation and corresponding $f_{ij}$ that satisfy this condition are given by \eqref{th2.4} and \eqref{th2.4a} with $\psi = 0$, i.e., given by Theorem \ref{teo7.4}. In that case, 
\begin{eqnarray*}
\Delta_{13} =  \pm \mu h(m_2+\lambda m_1 z_0),\quad \Delta_{23} = \mp  h(m_2+\lambda m_1 z_0),
\end{eqnarray*}
and then, by \eqref{3.15a}, $a=\pm 2/\mu$. Therefore, observing that $\Delta_{13}=-\mu \Delta_{23}$, \eqref{3.16a} reduces to 
$$
0=\mp 2\Delta_{13}+a\Delta_{23} = \mp 2(-\mu\Delta_{23})+(\pm 2/\mu)\Delta_{23}=\pm \frac{2(1+\mu^2)}{\mu}\Delta_{23}
$$ 
which holds if, and only if, $\Delta_{23}=0$, which implies $\Delta_{13}=0$ and, thus, a contradiction by \eqref{3.10a}. Hence, $c\neq 0$ on any open set.

From now on, we are assuming $c\neq 0$. If $a = 0$ on a open set, then \eqref{gauss} implies $b = \pm 1$ and \eqref{3.8} and \eqref{3.9} are equivalent to
\begin{eqnarray}
\mp 2\Delta_{13}-c\Delta_{23} = 0,&&\label{3.8b}\\
f_{21}D_t c -f_{22}D_x c -c\Delta_{13}\pm 2\Delta_{23} = 0,&&\label{3.9b}
\end{eqnarray}
It follows from \eqref{3.8b} and \eqref{3.10a} that $\Delta_{23}\neq 0$ and $c = \mp 2\Delta_{13}/ \Delta_{23}$. Since $\Delta_{13}$ and $\Delta_{23}$ depend only on $z_0$, $z_1$ and $z_2$, we conclude that $c$ depends only on $z_0$, $z_1$ and $z_2$ and \eqref{3.9b} reduces to
\begin{eqnarray}
f_{21}\sum_{i=0}^1 c_{z_i}z_{i,t}+f_{21}c_{z_2}[z_{0,t}-(\lambda z_0z_3+G)] -f_{22}\sum_{i=0}^2 c_{z_i}z_{i+1}-c\Delta_{13}\pm 2\Delta_{23} = 0.\label{3.9c}
\end{eqnarray}

Using \eqref{1} and taking the derivative of the latter expression with respect to $z_{0,t}$, $z_{1,t}$ and $z_{3}$ implies that
\begin{eqnarray}\label{3.9d}
f_{21} (c_{z_0}+c_{z_2}) = f_{21} c_{z_1} = (\lambda z_0f_{21} + f_{22})c_{z_2} = 0.
\end{eqnarray}
Replacing \eqref{3.9d} into \eqref{3.9c} we obtain
\begin{eqnarray}
-f_{21}c_{z_2}G -f_{22}\sum_{i=0}^1 c_{z_i}z_{i+1}-c\Delta_{13}\pm 2\Delta_{23} = 0.\label{3.9e}
\end{eqnarray}
Looking at \eqref{3.9d} is easy to see that, if $f_{21}\neq 0$ and $\lambda z_0f_{21}+f_{22}\neq 0$ then $c$ is a constant and from \eqref{3.8b} and \eqref{3.9b} we obtain $\Delta_{13}=\Delta_{23}=0$, which is a contradiction with \eqref{3.10a}. If $f_{21}= 0$ and $\lambda z_0f_{21}+f_{22}= 0$ on a open set, then, by \eqref{3.8b} and \eqref{3.9b}, $\Delta_{23} = 0$, which is a contradiction. Hence, we have only two possibilities, namely,\\
$$
(i)\hspace*{0.2 cm} f_{21} = 0 \quad \textnormal{and}\quad \lambda z_0f_{21}+f_{22}\neq 0,\qquad
(ii)\hspace*{0.2 cm} f_{21} \neq 0 \quad \textnormal{and}\quad \lambda z_0f_{21}+f_{22}= 0,
$$
on a open set. 

Assuming $(i)$, we can observe from \eqref{3.9d} and $\Delta_{12}\neq 0$ that $c_{z_2} = 0$ on a open set. Thus $f_{ij}$ are given by \eqref{th2.4a} with $\mu=m_1 = 0$ and $m_2\neq 0$ or \eqref{th2.5b} with $\mu = m_2 = 0$, i.e., given by Theorems \ref{teo7.4} with $\mu=m_1 = 0$ and $m_2\neq 0$ or \ref{teo7.5}-$(ii)$ with $\mu = m_2 = 0$, respectively.

If $f_{ij}$ are given by \eqref{th2.4a} with $\mu=m_1 = 0$ and $m_2\neq 0$ then $\Delta_{13} = 0$. However, $\Delta_{13}=0$ implies $c=0$ which is a contradiction, because we firstly showed that $c\neq 0$. If $f_{ij}$ are given by \eqref{th2.5b} with $\mu = m_2 = 0$ then
\begin{eqnarray}\label{3.9f}
\mp 2\left[ \pm (\lambda -\theta me^{\theta z_0})z_2 \mp m\theta^2e^{\theta z_0}z_1^2\right]+c\left(h+\frac{\theta}{a}\right)\left(\lambda -\theta me^{\theta z_0}\right)z_1 = 0.
\end{eqnarray}
Since $c_{z_2} = 0$, differentiating the latter equation with respect $z_2$ and, in following, with respect to $z_1$, and observing that $\lambda -\theta me^{\theta z_0}\neq 0$, we have $c h' = 0$ on a open set, which leads to a contradiction with \eqref{2.3.1}. This concludes $(i)$.

Assuming $(ii)$, i.e., $\lambda z_0f_{21}+f_{22}=\phi_{22}= 0$, we necessarily have $\phi_{12}\neq 0$, because $\Delta_{12} = -\phi_{12}f_{21}$. Moreover, from \eqref{3.9d} we conclude $c_{z_0}+c_{z_2}=c_{z_1}=0$.  Thus, $f_{ij}$ are given by \eqref{th2.2a} with $\mu=0$ or \eqref{th2.3a} with $\mu=0$ or \eqref{th2.4a} with $\mu \neq 0$, i.e., given by Theorems \ref{teo7.2} with $\mu=0$ or \eqref{teo7.3} with $\mu=0$ or \eqref{teo7.4} with $\mu \neq 0$, respectively.

If $f_{ij}$ are given by \eqref{th2.2a} with $\mu=0$ is easy to see that $\Delta_{13}=0$, which is a contradiction, because $c = \mp 2\Delta_{13}/ \Delta_{23}$ is not zero. If $f_{ij}$ are given by \eqref{th2.3a} with $\mu=0$ is easy to see that $\Delta_{13}=-\lambda \eta z_1(\neq 0)$ and $\Delta_{23} = -\lambda m_1 z_1(\neq 0)$. Therefore, it follows from $c = \mp 2\Delta_{13}/ \Delta_{23}$ that $c = \mp 2\eta /m_1$ is a constant and, from \eqref{3.8b} and \eqref{3.9b}, we obtain $\Delta_{13}=\Delta_{23}=0$, which is a contradiction with \eqref{3.10a}.

If $f_{ij}$ are given by \eqref{th2.4a} with $\mu\neq 0$ then $c = \pm 2\mu$  is a constant and, from \eqref{3.8b} and \eqref{3.9b}, we obtain $\Delta_{13}=\Delta_{23}=0$, which is a contradiction with \eqref{3.10a}. This concludes $(ii)$. Therefore, $a\neq 0$ on any open set and, thus, we conclude the proof of the Lemma \ref{lemma3.1}.

\cqd


Now, suppose that we have substituted the expressions of the total derivatives with respect to $x$ and $t$ given by \eqref{3.11} and \eqref{3.12} into equations \eqref{3.8} and \eqref{3.9}, i.e., 
\begin{eqnarray}\label{10}
f_{11}a_t + f_{21}b_t - f_{12}a_x - f_{22}b_x - 2b(f_{11}f_{32}-f_{12}f_{31})+(a-c)(f_{21}f_{32}-f_{22}f_{31})\hspace*{3 cm}\nonumber\\
-\sum_{i=0}^l (f_{12}a_{z_i}+f_{22}b_{z_i})z_{i+1}+\sum_{i=2}^l (f_{11}a_{z_i}+f_{21}b_{z_i})\partial_x^{i-2}(z_{0,t}-F) - \sum_{j=1}^m (f_{12}a_{w_j}+f_{22}b_{w_j})v_j\hspace*{1.2 cm}\nonumber\\
+\sum_{j=0}^m (f_{11}a_{w_j}+f_{21}b_{w_j})w_{j+1} - \sum_{k=1}^n (f_{12}a_{v_k}+f_{22}b_{v_k})\partial_t^{k-1}(z_{0,t}-F)+\sum_{k=0}^n (f_{11}a_{v_k}+f_{21}b_{v_k})v_{k+1}=0,
\end{eqnarray}
and
\begin{eqnarray}\label{11}
f_{11}b_t + f_{21}c_t - f_{12}b_x - f_{22}c_x + 2b(f_{21}f_{32}-f_{22}f_{31})+(a-c)(f_{11}f_{32}-f_{12}f_{31})\hspace*{3 cm}\nonumber\\
-\sum_{i=0}^l (f_{12}b_{z_i}+f_{22}c_{z_i})z_{i+1}+\sum_{i=2}^l (f_{11}b_{z_i}+f_{21}c_{z_i})\partial_x^{i-2}(z_{0,t}-F) - \sum_{j=1}^m (f_{12}b_{w_j}+f_{22}c_{w_j})v_j\hspace*{1.2 cm}\nonumber\\
+\sum_{j=0}^m (f_{11}b_{w_j}+f_{21}c_{w_j})w_{j+1} - \sum_{k=1}^n (f_{12}b_{v_k}+f_{22}c_{v_k})\partial_t^{k-1}(z_{0,t}-F)+\sum_{k=0}^n (f_{11}b_{v_k}+f_{21}c_{v_k})v_{k+1}=0.
\end{eqnarray}

If $m=n$, then differentiating \eqref{10} and \eqref{11} with respect to $v_{n +1}$ and $w_{n+1}$ leads to
\begin{eqnarray}\label{13}
\begin{array}{rr}
f_{11}a_{v_n}+f_{21}b_{v_n} = 0,\\ f_{11}b_{v_n}+f_{21}c_{v_n} = 0,
\end{array}
\begin{array}{rr}
f_{11}a_{w_n}+f_{21}b_{w_n} = 0,\\ f_{11}b_{w_n}+f_{21}c_{w_n} = 0.
\end{array}
\end{eqnarray}

If $f_{21}\neq 0$ on a non-empty open set (which is the case of the equations and $f_{ij}$ given by the Theorems \ref{teo7.2}, \ref{teo7.3} and \ref{teo7.5}-(i) and also may be the case of the equations given by Theorems \ref{teo7.4} and \ref{teo7.5}-(ii)) then

\begin{eqnarray}\label{13,a}
\begin{array}{rcl}
&& b_{v_n} = -\frac{f_{11}}{f_{21}}a_{v_n},\\ && c_{v_n} = \left(\frac{f_{11}}{f_{21}}\right)^2 a_{v_n},
\end{array}
\begin{array}{rcl}
&& b_{w_n} = -\frac{f_{11}}{f_{21}}a_{w_n},\\ && c_{w_n} = \left(\frac{f_{11}}{f_{21}}\right)^2 a_{w_n}.
\end{array}
\end{eqnarray}

Differentiating the Gauss equation with respect to $v_n$ and $w_n$ leads to $a_{v_n}c+ac_{v_n}-2bb_{v_n} = 0$ and $a_{w_n}c+ac_{w_n}-2bb_{w_n} = 0$, respectively, and using \eqref{13,a} in such derivatives we obtain
\begin{eqnarray}\label{1.18}
\left[ c + \left(\frac{f_{11}}{f_{21}}\right)^2 a+2\frac{f_{11}}{f_{21}} b\right] a_{v_n}=0, \quad \left[ c + \left(\frac{f_{11}}{f_{21}}\right)^2 a+2\frac{f_{11}}{f_{21}} b\right] a_{w_n}=0.
\end{eqnarray}

The equation \eqref{1.18} holds when $m=n$ and $f_{21}\neq 0$. The cases $m < n$ or $m > n$ need to be considered separately, and they will be analyzed in Lemmas \ref{lemma3.2}-\ref{lemma3.4}. The case $f_{21}\equiv 0$ will be considered in Lemma \ref{lemma3.4}. 

The discussing leading to \eqref{1.18} shows that the analysis of the Codazzi equations (\eqref{3.8} and \eqref{3.9}) splits naturally into several branches which are characterized by the vanishing or non-vanishing of $f_{21}$ and the expression between brackets in \eqref{1.18}. The various cases are treated in Lemmas \ref{lemma3.3}-\ref{lemma3.4} and are organized according to the figure below.

\begin{figure}[ht]
\begin{center}
\includegraphics[scale=0.5]{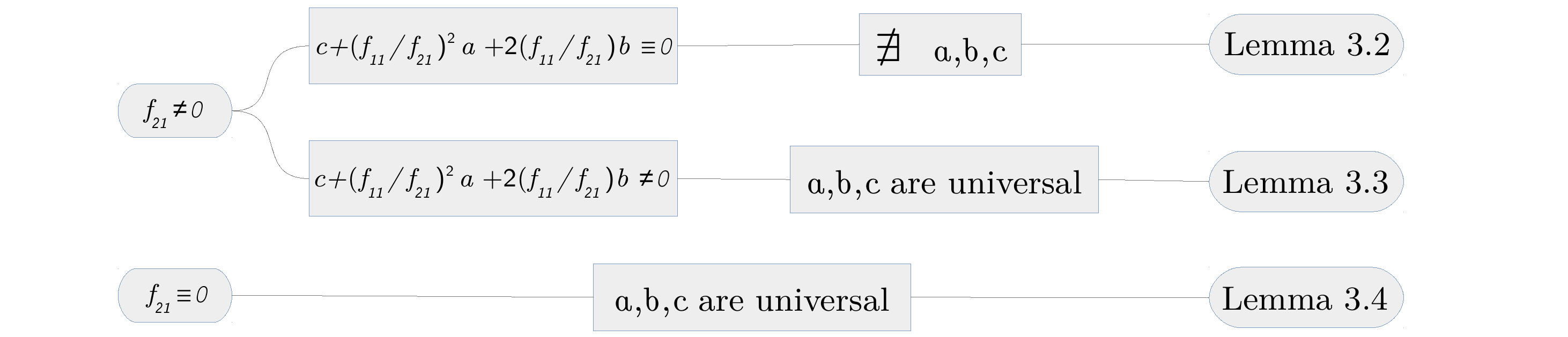}
\end{center}
\end{figure}


\begin{lemma}\label{lemma3.3}
Consider an equation of type \eqref{T} describing pseudospherical surfaces, under the condition \eqref{1}, as given by the Theorems $\ref{teo7.2}$-$\ref{teo7.5}$. Assume there is a local isometric immersion of the pseudospherical surface determined by a solution $u(x,t)$ of \eqref{T}, for which the coefficients $a$, $b$ and $c$ of the second fundamental form depend on $x, t, z_0, \ldots, z_l, w_1, \ldots, w_m, v_1, \ldots, v_n$, where $1\leq l  <\infty$, $1\leq m < \infty$ and $1\leq n < \infty$ are finite, but otherwise arbitrary. Suppose $f_{21}\neq 0$ on a non-empty open set. If
\begin{eqnarray}\label{19}
c + \left(\frac{f_{11}}{f_{21}}\right)^2 a+2\frac{f_{11}}{f_{21}} b = 0,
\end{eqnarray}
on a non-empty open set, then the equations \eqref{gauss}, \eqref{3.8} and \eqref{3.9} form an inconsistent system.
\end{lemma}

\noindent\textbf{Proof}. Firstly, let use \eqref{19} and the Gauss equation in order to obtain $b$ and $c$ in terms of $a$, $f_{11}$ and $f_{21}$. We will then substitute the total derivatives of $b$ and $c$ back into \eqref{3.8} and \eqref{3.9}. 

If \eqref{19} holds then substituting $c$ into the Gauss equation leads to
\begin{eqnarray}
&& b = \pm 1 - \frac{f_{11}}{f_{21}}a,\label{70}\\
&& c = \left(\frac{f_{11}}{f_{21}}\right)^2 a \mp 2\frac{f_{11}}{f_{21}}.\label{70.1}
\end{eqnarray}
Moreover, using \eqref{1} and \eqref{2.3.1} we can see that $(f_{11}/f_{21})_{,z_0}+(f_{11}/f_{21})_{,z_2} = 0$, and thus
\begin{eqnarray*}
f_{11}D_t a + f_{21}D_t b &=& f_{21} (\lambda z_0z_3+G)\left(\frac{f_{11}}{f_{21}}\right)_{,z_2}a,\\
f_{12}D_x a+f_{22}D_x b &=& -\frac{\Delta_{12}}{f_{21}}D_x a - f_{22}a \left[\left(\frac{f_{11}}{f_{21}}\right)_{,z_0}z_1 +\left(\frac{f_{11}}{f_{21}}\right)_{,z_2}z_3\right],\\
f_{11}D_t b+f_{21}D_t c &=& -(f_{11}a \mp 2f_{21})(\lambda z_0z_3+G)\left(\frac{f_{11}}{f_{21}}\right)_{,z_2},\\
f_{12}D_xb + f_{22}D_x c &=& \frac{f_{11}}{f_{21}}\frac{\Delta_{12}}{f_{21}}D_x a +\left[\frac{\Delta_{12}}{f_{21}}a+\frac{f_{22}}{f_{21}}(f_{11}a \mp 2f_{21})\right]\left[\left(\frac{f_{11}}{f_{21}}\right)_{,z_0}z_1 + \left(\frac{f_{11}}{f_{21}}\right)_{,z_2}z_3\right].                   
\end{eqnarray*}
Equation \eqref{3.8} becomes,
\begin{eqnarray}\label{5.41}
a[f_{21}G+(\lambda z_0f_{21}+f_{22})z_3]\left(\frac{f_{11}}{f_{21}}\right)_{,z_2}+af_{22}\left(\frac{f_{11}}{f_{21}}\right)_{,z_0}z_1 + \frac{\Delta_{12}}{f_{21}}D_xa -2b\Delta_{13}+(a-c)\Delta_{23} = 0,
\end{eqnarray}
and \eqref{3.9} becomes
\begin{eqnarray}\label{5.42}
-\frac{1}{f_{21}}(f_{11}a\mp 2 f_{21})[f_{21}G+(\lambda z_0 f_{21}+f_{22})z_3] \left(\frac{f_{11}}{f_{21}}\right)_{,z_2} -\frac{\Delta_{12}}{f_{21}}a\left[\left(\frac{f_{11}}{f_{21}}\right)_{,z_0}z_1+\left(\frac{f_{11}}{f_{21}}\right)_{,z_2}z_3\right]\hspace*{1 cm}\nonumber  \\
-\frac{f_{22}}{f_{21}}(f_{11}a\mp 2f_{21}) \left(\frac{f_{11}}{f_{21}}\right)_{,z_0}z_1 -\frac{f_{11}}{f_{21}}\frac{\Delta_{12}}{f_{21}}D_xa  +(a-c)\Delta_{13}+2b\Delta_{23} = 0.
\end{eqnarray}

Observe that in \eqref{5.41} and \eqref{5.42} we only have the total derivative of the coefficient $a$ of the second fundamental form of the local isometric immersion. We are going to use expressions \eqref{5.41} and \eqref{5.42} in an equivalent form that will be more convenient to work with when the expressions of $G$ and $f_{ij}$ given in Theorem \ref{teo7.2}-\ref{teo7.5} are taken into account.

Adding equation \eqref{5.41} multiplied by $f_{11}/f_{21}$ with \eqref{5.42}, we get
\begin{eqnarray}\label{5.43}
\left[ \pm 2f_{21}G\pm 2(f_{22}+\lambda z_0f_{21})z_3 -\frac{\Delta_{12}}{f_{21}}az_3\right]\left(\frac{f_{11}}{f_{21}}\right)_{,z_2} + \left[-\frac{\Delta_{12}}{f_{21}}a\pm 2 f_{22}\right]\left(\frac{f_{11}}{f_{21}}\right)_{,z_0}z_1\hspace*{1 cm}\nonumber\\
+\left[1+\left(\frac{f_{11}}{f_{21}}\right)^2\right]\left[a\Delta_{13}+\left(-\frac{f_{11}}{f_{21}}a \pm 2\right)\Delta_{23} \right] = 0.
\end{eqnarray}
Taking the $v_k$ and $w_j$, $1\leq k\leq n$ and $1\leq j\leq m$, derivatives of \eqref{5.43}, we have, respectively,
\begin{eqnarray}\label{5.43b}
Q a_{v_k} = 0, \quad Q a_{w_j} = 0,
\end{eqnarray}
where
\begin{eqnarray}\label{5.43a}
Q = -\frac{\Delta_{12}}{f_{21}}\left(\frac{f_{11}}{f_{21}}\right)_{,z_2}z_3 -\frac{\Delta_{12}}{f_{21}}\left(\frac{f_{11}}{f_{21}}\right)_{,z_0}z_1+\left[1+\left(\frac{f_{11}}{f_{21}}\right)^2\right]\left(\Delta_{13}-\frac{f_{11}}{f_{21}}\Delta_{23}  \right).
\end{eqnarray}

Suppose $Q = 0$ on a non-empty open set. Differentiating $Q$ with respect to $z_3$ we have $(f_{11}/f_{21})_{,z_2} = 0$ and, consequently, $(f_{11}/f_{21})_{,z_0} = 0$. Hence $f_{11}/f_{21}$ is a nonzero constant, which happens only in the branches of the classification corresponding to Theorems \ref{teo7.4} with $m_1=0$ and \ref{teo7.5}-$(ii)$ with $\eta = 0$.

If $f_{ij}$ are given by \eqref{th2.4a}, i.e., Theorem \ref{teo7.4}, with $m_1 = 0$, then
$$
\frac{f_{11}}{f_{21}}=\frac{1}{\mu}, \quad \Delta_{12}=m_2\sqrt{1+\mu^2}, \quad \Delta_{13}=\pm m_2 \mu h, \quad \Delta_{23} = \mp m_2h.
$$
Therefore, \eqref{5.43a} implies that $\Delta_{13}-f_{11}\Delta_{23}/f_{21}=0$ if, and only if, $m_2 = 0$ and, thus, a contradiction.

On the other hand, if $f_{ij}$ are given by \eqref{th2.5b}, i.e., Theorem \ref{teo7.5}-$(ii)$, with $\eta = 0$, then
$$
\frac{f_{11}}{f_{21}}=\frac{1}{\mu}, \quad \Delta_{12}=\pm \sqrt{1+\mu^2}f_{11}z_1, \quad \Delta_{13}=(m_1\theta e^{\theta z_0}-\lambda)f_{11}\left(-\mu z_1\pm \frac{1}{A\sqrt{1+\mu^2}}\right)\mp \frac{\theta}{A\sqrt{1+\mu^2}},
$$
$$
\Delta_{23} = (m_1\theta e^{\theta z_0}-\lambda)f_{11}\left( z_1 \pm \frac{\mu}{A\sqrt{1+\mu^2}}\right)\mp \frac{\theta}{A\sqrt{1+\mu^2}}\phi_{22}.
$$
Therefore, \eqref{5.43a} is equivalent to
\begin{eqnarray}\label{(54)}
0=\mu \Delta_{13}-\Delta_{23} = -(m_1\theta e^{\theta z_0}-\lambda)(1+\mu^2)z_1f_{11}\pm \frac{\theta}{A\sqrt{1+\mu^2}}(\phi_{22}-\mu),
\end{eqnarray}
where by \eqref{(7.2)} we know that
$$
\phi_{22} = f_{22}+\lambda z_0f_{21}= \mu a  m_1\theta e^{\theta z_0}z_1^2
+ (m_1\theta e^{\theta z_0}-\lambda) \frac{1}{\theta}
\left[\mu(az_0+b) \mp \frac{\theta z_1}{\sqrt{1+\mu^2}}
  \right].
$$
Differentiating \eqref{(54)} with respect to $z_2$ leads to $-(m_1\theta e^{\theta z_0}-\lambda)(1+\mu^2)z_1f_{11,z_2} = 0$, which holds if, and only if, $m_1\theta e^{\theta z_0}-\lambda=0$, i.e., if, and only if, $\lambda=m_1 = 0$ and, thus, a contradiction. Hence, we have shown that $Q$ does not vanish in a non-empty open set.

Let $Q\neq 0$, on a non-empty open set. We are going to showing that this also leads us into some contradiction. Consequently, \eqref{5.43b} implies that $a_{v_k} = a_{w_j} = 0$, $k =1,2, \ldots,n$ and $j=1,2, \ldots, m$ and, thus, $a$ is a function depending only on $x, t, z_0, \ldots, z_l$. But, differentiating \eqref{5.43} with respect to $z_l$, $l\geq 4$, we also obtain $Qa_{z_l} = 0$ where $Q$ is given by \eqref{5.43a} and, since $Q\neq 0$, we conclude that $a$ depends only on $x, t, z_0, \ldots, z_3$. Moreover, differentiating \eqref{5.41} with respect to $z_4$ leads to
\begin{eqnarray*}
\frac{\Delta_{12}}{f_{21}}a_{z_3} = 0,
\end{eqnarray*}
i.e., $a_{z_3} = 0$, on a open set.

Taking the $z_3$ derivative of \eqref{5.43}, we obtain
\begin{eqnarray}\label{5.44}
\left[-\frac{\Delta_{12}}{f_{21}}a\pm 2(f_{22}+\lambda z_0f_{21})\right]\left(\frac{f_{11}}{f_{21}}\right)_{,z_2}=0.
\end{eqnarray}

Suppose $\left(f_{11}/f_{21}\right)_{,z_2} \equiv 0$, which happens only in \eqref{th2.4a} with $m_1=0$ or \eqref{th2.5b} with $\eta = 0$, i.e., in the branches of the classification corresponding to Theorems \ref{teo7.4} with $m_1=0$ and \ref{teo7.5}-$(ii)$ with $\eta=0$. Therefore, if $f_{ij}$ are given by \eqref{th2.4a} with $m_1 = 0$ then
$$
\frac{f_{11}}{f_{21}}=\frac{1}{\mu}, \quad \Delta_{12}=m_2\sqrt{1+\mu^2}h, \quad \Delta_{13}=\pm m_2 \mu h, \quad \Delta_{23} = \mp m_2h.
$$
Thus, \eqref{5.44} is satisfied and, substituting back into \eqref{5.43}, we get
\begin{eqnarray*}
0=a\Delta_{13}+\left(-\frac{f_{11}}{f_{21}}a \pm 2\right)\Delta_{23} = a\left(\frac{1+\mu^2}{\mu}\right) \mp2,
\end{eqnarray*}
i.e., $a$ is a constant. But if $a$ is a contant, it follows from \eqref{5.41} and \eqref{5.42} that
\begin{eqnarray*}
-2b\Delta_{13}+(a-c)\Delta_{23} = 0,&&\\
(a-c)\Delta_{13}+2b\Delta_{23} = 0,&&
\end{eqnarray*}
which implies that $\Delta_{13}=\Delta_{23}=0$ and, thus, a contradiction with \eqref{3.10a}. 

On the other hand, considering the case in which $f_{ij}$ are given by \eqref{th2.5b} with $\eta = 0$ we get
$$
\frac{f_{11}}{f_{21}}=\frac{1}{\mu}, \quad \Delta_{12}=\pm \sqrt{1+\mu^2}f_{11}z_1, \quad \Delta_{13}=(m_1\theta e^{\theta z_0}-\lambda)f_{11}\left(-\mu z_1\pm \frac{1}{A\sqrt{1+\mu^2}}\right)\mp \frac{\theta}{A\sqrt{1+\mu^2}},
$$
$$
\Delta_{23} = (m_1\theta e^{\theta z_0}-\lambda)f_{11}\left( z_1 \pm \frac{\mu}{A\sqrt{1+\mu^2}}\right)\mp \frac{\theta}{A\sqrt{1+\mu^2}}\phi_{22}.
$$
Therefore, \eqref{5.44} is satisfied and, substituting back into \eqref{5.43}, we get
\begin{eqnarray}\label{5.45}
0=a\Delta_{13}+\left(-\frac{a}{\mu} \pm 2\right)\Delta_{23} =\frac{1}{\mu} \left[ a(\mu \Delta_{13}-\Delta_{23})\pm 2\mu \Delta_{23}\right],
\end{eqnarray}
and since 
\begin{eqnarray*}
0\neq \mu \Delta_{13}-\Delta_{23} = -(m_1\theta e^{\theta z_0}-\lambda)(1+\mu^2)z_1f_{11}\pm \frac{\theta}{A\sqrt{1+\mu^2}}(\phi_{22}-\mu),
\end{eqnarray*}
on a open set, we observe that $a$ depends only on $z_0$, $z_1$ and $z_2$. Then, from \eqref{5.41} and \eqref{5.42}, we obtain
\begin{eqnarray}
\frac{\Delta_{12}}{f_{21}}D_xa -2b\Delta_{13}+(a-c)\Delta_{23} = 0,&&\label{5.46}\\
-\frac{1}{\mu}\frac{\Delta_{12}}{f_{21}}D_xa  +(a-c)\Delta_{13}+2b\Delta_{23} = 0.&&\label{5.47}
\end{eqnarray}

Differentiating \eqref{5.46} with respect to $z_3$ leads to $\Delta_{12}a_{z_2}/f_{21}=0$ and, therefore, $a_{z_2}=0$. Differentiating \eqref{5.45} with respect to $z_2$ and replacing the result back into \eqref{5.42} leads to
\begin{eqnarray}
-a(1+\mu^2)z_1 \pm 2\mu \left(z_1\pm \frac{\mu}{A\sqrt{1+\mu^2}}\right)=0,&&\label{5.48}\\
(\pm a-2\mu)\theta \phi_{22}\mp a\mu = 0.&&\label{5.49}
\end{eqnarray}
It follows from \eqref{5.46} and \eqref{5.47} that $a$ can not be constant, otherwise,
\begin{eqnarray*}
 -2b\Delta_{13}+(a-c)\Delta_{23} = 0,&&\\
(a-c)\Delta_{13}+2b\Delta_{23} = 0,&&
\end{eqnarray*}
and thus $\Delta_{13}=\Delta_{23}=0$, contradicting \eqref{3.10a}. Hence, \eqref{5.48} and \eqref{5.49} imply that $\phi_{22}$ depends only on $z_1$, i.e.,
\begin{eqnarray*}
0=\phi_{22,z_0}= \mu Am_1\theta^2 e^{\theta z_0}z_1^2+m_1\theta e^{\theta z_0} \left[\mu(Az_0+b)\mp \frac{\theta z_1}{\sqrt{1+\mu^2}}\right]+(m_1\theta e^{\theta z_0}-\lambda)\frac{\mu A}{\theta}.
\end{eqnarray*}
Differentiating the latter expression twice in $z_1$ we obtain $\mu Am_1\theta^2 e^{\theta z_0}=0$ which holds if, and only if, $m_1=0$. But, $m_1=0$ in the latter equation gives us $-\lambda \mu A/\theta=0$ which implies $\lambda = 0$. Hence, $m_1=\lambda = 0$ which is a contradiction with $\Delta_{12}\neq 0$. Thus, we have shown that $\left(f_{11}/f_{21}\right)_{,z_2}$ does not vanish on a non-empty open set.

Let us go back to equation \eqref{5.44} and analyze the case $\left(f_{11}/f_{21}\right)_{,z_2} \neq 0$ on a non-empty open set. This condition holds in \eqref{th2.2a} or \eqref{th2.3a} or \eqref{th2.4a} with $m_1\neq 0$ or \eqref{th2.5a} or \eqref{th2.5b} with $\eta\neq 0$, i.e., in the branches of the classification corresponding to Theorems \ref{teo7.2}-\ref{teo7.5}. 

From \eqref{5.44}, $a$ is given by
\begin{eqnarray}\label{5.50}
a = \pm 2\frac{\phi_{22}f_{21}}{\Delta_{12}},
\end{eqnarray}
where, by \eqref{(7.2)}, we know that  $\phi_{22}=f_{22}+\lambda z_0f_{21}$. Thus, usando \eqref{5.50} and the fact of $\Delta_{13}-f_{11}\Delta_{23}/f_{21}=f_{31}\Delta_{12}/f_{21}$, the equations \eqref{5.41} and \eqref{5.43} are equivalent to
\begin{eqnarray}\label{5.51}
 a(f_{21}G  -f_{22}z_1)\left(\frac{f_{11}}{f_{21}}\right)_{,z_2} + \frac{\Delta_{12}}{f_{21}}\sum_{i=0}^1 a_{z_i}z_{i+1}+a \left[1+\left(\frac{f_{11}}{f_{21}}\right)^2\right]\Delta_{23}-2\left[\pm 1-a\frac{f_{11}}{f_{21}}\right]\frac{f_{31}}{f_{21}}\Delta_{12} = 0,
\end{eqnarray}
\begin{eqnarray}\label{5.52}
 \pm 2f_{21}G\left(\frac{f_{11}}{f_{21}}\right)_{,z_2} \mp 2 \lambda z_0z_1f_{21}\left(\frac{f_{11}}{f_{21}}\right)_{,z_0}+\left[1+\left(\frac{f_{11}}{f_{21}}\right)^2\right]\left(a\frac{f_{31}}{f_{21}}\Delta_{12}\pm 2\Delta_{23} \right) = 0.
\end{eqnarray}
Remember that $\left(f_{11}/f_{21}\right)_{,z_0}+\left(f_{11}/f_{21}\right)_{,z_2} = 0$ and, by \eqref{(7.2)}, we also have $af_{31}\Delta_{12}/f_{21}\pm 2\Delta_{23} = \pm 2f_{21}\phi_{32}$. Hence, it follows from \eqref{5.52} that
\begin{eqnarray}\label{5.53}
G = -\lambda z_0z_1 - \frac{(f_{11}^2+f_{21}^2)}{f_{21}^2}\frac{\phi_{32}}{L},\quad \textnormal{where}\quad L:=\left(\frac{f_{11}}{f_{21}}\right)_{,z_2}.
\end{eqnarray}

If $G$ and $f_{ij}$ are given as in \eqref{th2.2} and \eqref{th2.2a}, i.e., Theorem \ref{teo7.2}, from \eqref{5.53} we get $L = -m\sqrt{1+\mu^2}h'/f_{21}^2$ and
\begin{eqnarray}\label{5.54}
 z_1\psi_{,z_0} + z_2 \psi_{,z_1} \pm m\psi =G=  \pm\frac{1}{m}(f_{11}^2+f_{21}^2)\psi.
\end{eqnarray}
Differentiating the latter with respect to $z_2$, there exists a function $P=P(z_0)$ such that
\begin{eqnarray*}\label{5.55}
 \frac{\psi_{,z_1}}{\psi} =  \pm\frac{1}{m}(f_{11}^2+f_{21}^2)_{,z_2}=P,
\end{eqnarray*}
i.e., there exist functions $R=R(z_0)$ and $S=S(z_0)$ such that
\begin{eqnarray}
\psi &=& Re^{Pz_1},\quad R\neq 0,\label{5.55}\\
\pm \frac{1}{m}(f_{11}^2+f_{21}^2) &=& Pz_2 + S.\label{5.56}
\end{eqnarray}
Differentiating \eqref{5.56} with respect $z_0$ and adding the result with the $z_2$ derivative of \eqref{5.56}, and using $f_{11,z_0}+f_{11,z_2}=0$ and $f_{21,z_0}+f_{21,z_2}=0$, we obtaing $P=A$ and $S = -Az_0+C$, where $A$ and $C$ are constants with $A\neq 0$. In fact, if $A=0$ then $P=0$ and $S=C$, and differentiating \eqref{5.56} with respect to $z_2$ leads to
$$
-(1+\mu^2)h \mp m\mu\sqrt{1+\mu^2} = 0,
$$
because $h'\neq 0$ on a open set. But, differentiating the latter again with respect to $z_2$ we get $(1+\mu^2)h'=0$ and, thus, a contradiction. So, $A\neq 0$.

Substituting \eqref{5.55} and \eqref{5.56} into \eqref{5.54}, we have
\begin{eqnarray}\label{l}
z_1 R' +(\pm m +Az_0 -C)R = 0.
\end{eqnarray}
Differentiating the latter with respect $z_1$ we get $R=constant$ which, when replaced back in \eqref{l}, gives us $AR=0$, which implies $R=0$ and, thus, a contradiction.

If $G$ and $f_{ij}$ are given as Theorem \ref{teo7.3}, then from \eqref{5.53} we obtain $L = -\eta h'/f_{21}^2$ and
$$
\eta^2( z_1h + m_1z_1 +m_2z_2) = (f_{11}^2+f_{21}^2)[m_1(1+\mu^2)-\mu\eta]z_1
$$
Taking the $z_1$ derivative of the last equation and replacing the result back into the same equation we obtain $\eta^2 m_2 = 0$, which contradicts the condition $\eta m_2\neq 0$ appearing in Theorem \ref{teo7.3}.

If $G$ and $f_{ij}$ are given as Theorem \ref{teo7.4} with $m_1\neq 0$, then from \eqref{5.53} we obtain $L = -m_1\sqrt{1+\mu^2}h'/f_{21}^2$ and
\begin{eqnarray}\label{(68)a}
-(\lambda z_1\pm \lambda m_1 z_0\pm m_2)h+z_1\psi_{,z_0}+z_2\psi_{,z_1}\pm m_1\psi = \frac{f_{11}^2+f_{21}^2}{m_1\sqrt{1+\mu^2}}[\pm \sqrt{1+\mu^2}\psi \pm \lambda m_1\mu z_0\pm m_2\mu].
\end{eqnarray}
Differentiating \eqref{(68)a} with respect to $z_0$ and $z_2$ and adding both results lead to
\begin{eqnarray}\label{(69)a}
\mp \lambda m_1h + z_1\psi_{,z_0z_0} + z_2 \psi_{,z_1z_0}\pm m_1\psi_{,z_0}+\psi_{,z_1} = \frac{f_{11}^2+f_{21}^2}{m_1\sqrt{1+\mu^2}}[\pm \sqrt{1+\mu^2}\psi_{,z_0} \pm \lambda m_1\mu ].
\end{eqnarray}
Likewise, differentiating \eqref{(69)a} with respect to $z_0$ and $z_2$ and adding both results, we obtain
\begin{eqnarray}\label{(70)a}
 z_1\psi_{,z_0z_0z_0} + z_2 \psi_{,z_1z_0z_0}\pm m_1\psi_{,z_0z_0}+2\psi_{,z_1z_0} = \pm\frac{f_{11}^2+f_{21}^2}{m_1}\psi_{,z_0z_0}.
\end{eqnarray}
Taking the $z_2$ derivative of \eqref{(70)a}, we have
\begin{eqnarray}\label{(71)a}
\psi_{,z_1z_0z_0} = \pm \frac{2}{m_1}(f_{11}f_{11,z_2}+f_{21}f_{21,z_2})\psi_{,z_0z_0}.
\end{eqnarray}
We now divide our analysis in two cases. According to whether $\psi_{,z_0z_0}\equiv 0$ or $\psi_{,z_0z_0}\neq 0$. 

If $\psi_{,z_0z_0}\equiv 0$ then we have from \eqref{(71)a} that $\psi_{,z_1z_0z_0}=0$ and, by \eqref{(70)a}, $\psi_{,z_1z_0} = 0$. Hence, $\psi = Az_0 + N$, where $A$ is a constant and $N=N(z_1)$ is a differentiable function. It follows from \eqref{(69)a} that
\begin{eqnarray}\label{(72)a}
\mp \lambda m_1h \pm m_1\psi_{,z_0}+\psi_{,z_1} = \frac{f_{11}^2+f_{21}^2}{m_1\sqrt{1+\mu^2}}[\pm \sqrt{1+\mu^2}A \pm \lambda m_1\mu ].
\end{eqnarray}
Differentiating \eqref{(72)a} with respect to $z_2$, since $h'\neq 0$, leads to
$$
\mp \lambda m_1 = 2\frac{\pm \sqrt{1+\mu^2}A \pm \lambda m_1\mu }{m_1\sqrt{1+\mu^2}}(-f_{11}-\mu f_{21}).
$$
Differentiating the last equation with respect to $z_2$ implies that $\pm \sqrt{1+\mu^2}A \pm \lambda m_1\mu =0$ and, from the last equation again, we have $\lambda m_1 = 0$. Since $m_1\neq 0$ we have $\lambda =0$. Thus, $A=0$. Hence, $\psi = N$. But, from \eqref{(69)a}, since $A=\lambda =0$, we have $N'=0$, i.e., $N$ is a constant. Finally, the equation \eqref{(68)a} gives us
\begin{eqnarray}\label{(73)a}
\mp m_2h \pm m_1N = \frac{f_{11}^2+f_{21}^2}{m_1\sqrt{1+\mu^2}}[\pm \sqrt{1+\mu^2}N \pm m_2\mu].
\end{eqnarray}
Taking the $z_2$ derivative of the last equation, since $h'\neq 0$, we get
$$
\mp m_2 = 2\frac{\pm \sqrt{1+\mu^2}N \pm  m_2\mu }{m_1\sqrt{1+\mu^2}}(-f_{11}-\mu f_{21}).
$$
Differentiating the latter equation with respect to $z_2$ leads to $\pm \sqrt{1+\mu^2}N \pm  m_2\mu=0$ and, thus, $m_2=0$. But, $\lambda = m_2 = 0$ contradicts the fact of $\Delta_{12}\neq 0$. Hence, we have shown that from equation \eqref{(71)a} we can not have $\psi_{,z_0z_0}\equiv 0$.

Let us now consider the case $\psi_{,z_0z_0}\neq 0$ in \eqref{(71)a}. So, it follows from \eqref{(71)a} that
\begin{eqnarray}\label{m}
\frac{\psi_{,z_1z_0z_0}}{\psi_{,z_0z_0}} = \pm \frac{2}{m_1}(f_{11}f_{11,z_2}+f_{21}f_{21,z_2}) = R(z_0),
\end{eqnarray}
where $R=R(z_0)$ is a differentiable function. Equation \eqref{m} may be written as
\begin{eqnarray}
\psi_{,z_1z_0z_0} &=& R\psi_{,z_0z_0},\label{(74,1)}\\
f_{11}^2+f_{21}^2 &=& \pm m_1 R z_2 + S(z_0),\label{(75,1)}
\end{eqnarray}
where $S=S(z_0)$ is a differentiable function. Taking the $z_0$ and $z_2$ derivative of \eqref{(75,1)}, adding the result and using $f_{11,z_0}+f_{11,z_2} = 0$ and $f_{21,z_0}+f_{21,z_2} = 0$ we obtain $R = -A$ constant and $S = \pm Am_1 z_0+B$ with $B$ constant. Hence,
\begin{eqnarray}
f_{11}^2+f_{21}^2 &=& \pm m_1 A (z_0- z_2) + B,\label{(76,1)}
\end{eqnarray}
and integrating once with respect to $z_0$ the equation \eqref{(74,1)}, we get
\begin{eqnarray}
\psi_{,z_1z_0} &=& -A\psi_{,z_0}+T(z_1),\label{(77,1)}
\end{eqnarray}
where $T=T(z_1)$ is a differentiable function. Substituting \eqref{(76,1)} and \eqref{(77,1)} into \eqref{(69)a} leads to
{\small \begin{eqnarray}\label{(78,1)}
\mp \lambda m_1h + z_1\psi_{,z_0z_0} + z_2[-A\psi_{,z_0}+T(z_1)]\pm m_1\psi_{,z_0}+\psi_{,z_1} = \frac{[\pm m_1 A (z_0- z_2) + B]}{m_1\sqrt{1+\mu^2}}[\pm \sqrt{1+\mu^2}\psi_{,z_0} \pm \lambda m_1\mu ].
\end{eqnarray}}
Taking the $z_2$ derivative of \eqref{(78,1)}, we have
\begin{eqnarray}\label{(79,1)}
\pm \lambda m_1h' = -T(z_1) \mp \frac{\lambda m_1\mu A}{\sqrt{1+\mu^2}} = \lambda m_1 C,
\end{eqnarray}
where $C$ is a nonzero constant, since $h'\neq 0$. Thus, from \eqref{(79,1)} we obtain $f_{11} = h = \pm C(z_0-z_2)+D$, where $D$ is a constant. But, replacing $f_{11}$ into \eqref{(76,1)} and using $f_{21} = \mu f_{11}\pm m_1\sqrt{1+\mu^2}$ and differentiating the remainder expression twice with respect to $z_0$ leads to $C=0$ and, thus, a contradiction with $h'\neq 0$. Therefore, we have shown that from equation \eqref{(71)a} we can not have $\psi_{,z_0z_0}\neq 0$ neither, on a non-empty open set. So, \eqref{5.53} is not true if $G$ and $f_{ij}$ are given as Theorem \ref{teo7.4} with $m_1\neq 0$.

Now, let consider $G$ and $f_{ij}$ given by Theorem \ref{teo7.5}-$(i)$. Then, from \eqref{5.53} we obtain $L = -A\eta/f_{21}^2$ and
\begin{eqnarray}\label{(74)a}
\lambda \left(z_1z_2 - 2 z_0z_1 - \frac{m}{\tau}z_1 \mp  \frac{z_2}{\tau}\right)
+ \tau e^{\pm \tau z_1}\left( \tau z_0z_2 \pm z_1+ m z_2\right)\varphi\hspace*{4 cm}\nonumber\\
\pm e^{\pm \tau z_1}\left(\tau z_0z_1 +\tau z_1z_2 + mz_1\pm z_2\right)\varphi' + z_1^2 e^{\pm \tau z_1}\varphi''=-\lambda z_0z_1+\frac{f_{11}^2+f_{21}^2}{A\eta}\phi_{32}.
\end{eqnarray}
Taking the $z_2$ derivative of \eqref{(74)a}, we obtain
\begin{eqnarray}\label{(74)a}
\lambda \left(z_1  \mp  \frac{1}{\tau}\right)
+ \tau e^{\pm \tau z_1}\left( \tau z_0 + m \right)\varphi \pm e^{\pm \tau z_1}\left(\tau z_1 \pm 1\right)\varphi' =-2\frac{f_{11}+\mu f_{21}}{\eta}\phi_{32},
\end{eqnarray}
which the derivative with respect to $z_2$ leads to $0=(1+\mu^2)A$, i.e., $A=0$, which contradicts $f_{11,z_2}\neq 0$.

Finally, if $G$ and $f_{ij}$ are given by Theorem \ref{teo7.5}-$(ii)$ with $\eta\neq 0$ then, from \eqref{5.53}, we have $L = -A\eta/f_{21}^2$ and
\begin{eqnarray}\label{(76)a1}
\lambda(2z_1z_2-3z_0z_1-m_2z_1) +
m_1\theta e^{\theta z_0}(\theta z_1^3+z_1z_2+2z_0z_1+m_2z_1) = -\lambda z_0z_1+\frac{f_{11}^2+f_{21}^2}{A\eta}\phi_{32},
\end{eqnarray}
where
$$
\phi_{32} = \pm \sqrt{1+\mu^2} A m_1\theta e^{\theta z_0}z_1^2
-(m_1\theta e^{\theta z_0}-\lambda) \frac{1}{\theta}
\left\{A\eta z_1 \mp\frac{1}{\sqrt{1+\mu^2}}
\left[(1+\mu^2)(Az_0+B)+\mu\eta +\frac{\theta}{A}\right]\right\} .
$$
Differentiating \eqref{(76)a1} three times with respect to $z_1$, we obtain $m_1\theta^2 e^{\theta z_0}=0$, i.e., $m_1=0$ (and then $\lambda\neq 0$). Thus, we can rewrite \eqref{(76)a1} such as
\begin{eqnarray}\label{(76)a2}
2z_1z_2-2z_0z_1-m_2z_1 = \frac{f_{11}^2+f_{21}^2}{A\theta\eta}
 \left\{A\eta z_1 \mp\frac{1}{\sqrt{1+\mu^2}}
\left[(1+\mu^2)(Az_0+B)+\mu\eta +\frac{\theta}{A}\right]\right\}.
\end{eqnarray}
Differentiating \eqref{(76)a2} with respect to $z_1$ leads to $f_{11}^2+f_{21}^2 = -2\theta (z_0-z_2)-m_2\theta$, which replaced back into \eqref{(76)a2} gives us
$$
(1+\mu^2)(Az_0+B)+\mu\eta +\frac{\theta}{A}=0.
$$
The $z_0$ derivative of the last equation implies that $(1+\mu^2)A = 0$, i.e., $A=0$, which contradicts $f_{11,z_2}\neq 0$. This concludes the proof of Lemma \ref{lemma3.3}.

\cqd


In the next two lemmas (Lemmas \ref{lemma3.2}-\ref{lemma3.4}) we will see that, under certain conditions, if a local isometric immersion exists for which the components $a$, $b$, $c$ of the second fundamental form depends only on a jet of finite order of $u$, then its coefficients are functions depending only on $x$ and $t$. Moreover, the proof in both lemmas requires separate the analysis of the cases $m=n$, $m < n$ and $n < m$.

\begin{lemma}\label{lemma3.2}
Consider an equation of type \eqref{T} describing pseudospherical surfaces, under the condition \eqref{1}, given by the Theorems $\ref{teo7.2}$-$\ref{teo7.5}$. Assume there is a local isometric immersion of the pseudospherical surface, determined by a solution $u(x,t)$ of \eqref{T}, for which the coefficients $a$, $b$ and $c$ of the second fundamental form depend on $x, t, z_0, \ldots, z_l, w_1, \ldots, w_m, v_1, \ldots, v_n$, where $1\leq l  <\infty$, $1\leq m < \infty$ and $1\leq n < \infty$ are finite, but otherwise arbitrary. Suppose $f_{21}\neq 0$ on a non-empty open set. If
\begin{eqnarray}\label{19.1}
c + \left(\frac{f_{11}}{f_{21}}\right)^2 a+2\frac{f_{11}}{f_{21}} b\neq 0
\end{eqnarray}
holds on a non-empty open set then $a$, $b$ and $c$ are functions of $x$ and $t$ only, and therefore universal.
\end{lemma}

\noindent\textbf{Proof}. Our analysis consists in three cases, namely,
$$
(i)\hspace{0.2 cm} m=n, \qquad (ii)\hspace{0.2 cm} m < n, \qquad (iii)\hspace{0.2 cm} n < m.
$$

Firstly, we consider the case $m=n$ and we are going to show that, from \eqref{10} and \eqref{11}, we have $a$, $b$ and $c$ depending only on $x$ and $t$.

Suppose $l =1$. If \eqref{19.1} holds then it follows from \eqref{1.18} that $a_{v_n}=0$ and $a_{w_n}=0$ and, consequently, by \eqref{13} we obtain $b_{v_n}=c_{v_n}=0$ and $b_{w_n}=c_{w_n}=0$. Thus, successive differentiation of \eqref{10}, \eqref{11} and \eqref{gauss} with respect to $v_{n+1}, \ldots, v_1$ and $w_{n+1}, \ldots, w_1$ lead to $a_{v_k}=b_{v_k}=c_{v_k}=0$ and $a_{w_k}=b_{w_k}=c_{w_k}=0$ for $k=0,1, \ldots, n$. Therefore, $a$, $b$ and $c$ are universal.

Suppose $l \geq 2$. Successive differentiation of \eqref{10}, \eqref{11} and \eqref{gauss} with respect to $v_{n+1}, \ldots, v_2$ and $w_{n+1}, \ldots, w_2$ lead to $a_{w_k}=b_{w_k}=c_{w_k}=0$ and $a_{v_k}=b_{v_k}=c_{v_k}=0$ for $k=1,2, \ldots, n$. In particular, $a$, $b$ and $c$ do not depend on $w_k$ and neither $v_k$ for $k=1,2, \ldots, n$. Therefore, $a$, $b$ and $c$ are functions of $x, t, z_0=w_0, z_1=v_0, \ldots, z_l$. Moreover, the equations \eqref{10} and \eqref{11} are equivalent to
\begin{eqnarray}\label{20}
f_{11}a_t + f_{21}b_t - f_{12}a_x - f_{22}b_x - 2b(f_{11}f_{32}-f_{12}f_{31})+(a-c)(f_{21}f_{32}-f_{22}f_{31})\hspace*{2 cm}\nonumber\\
-\sum_{i=0}^l (f_{12}a_{z_i}+f_{22}b_{z_i})z_{i+1}+\sum_{i=2}^l (f_{11}a_{z_i}+f_{21}b_{z_i})\partial_x^{i-2}(z_{0,t}-F)\hspace*{2 cm}\nonumber\\
+(f_{11}a_{w_0}+f_{21}b_{w_0})w_1 + (f_{11}a_{v_0}+f_{21}b_{v_0})v_1=0,
\end{eqnarray}
and
\begin{eqnarray}\label{21}
f_{11}b_t + f_{21}c_t - f_{12}b_x - f_{22}c_x + 2b(f_{21}f_{32}-f_{22}f_{31})+(a-c)(f_{11}f_{32}-f_{12}f_{31})\hspace*{2 cm}\nonumber\\
-\sum_{i=0}^l (f_{12}b_{z_i}+f_{22}c_{z_i})z_{i+1}+\sum_{i=2}^l (f_{11}b_{z_i}+f_{21}c_{z_i})\partial_x^{i-2}(z_{0,t}-F)\hspace*{2 cm}\nonumber\\
+(f_{11}b_{w_0}+f_{21}c_{w_0})w_1 + (f_{11}b_{v_0}+f_{21}c_{v_0})v_1=0.
\end{eqnarray}

Differentiating \eqref{20} and \eqref{21} with respect to $z_{\ell+1}$, we obtain, respectively,
\begin{eqnarray*}
(f_{12}+\lambda z_0f_{11})a_{z_\ell}+(f_{22}+\lambda z_0f_{21})b_{z_\ell}=0,\quad (f_{12}+\lambda z_0f_{11})b_{z_\ell}+(f_{22}+\lambda z_0f_{21})c_{z_\ell}=0,
\end{eqnarray*}
and, using \eqref{(7.2)}, we have
\begin{eqnarray}
\begin{array}{rr}\label{31}
\phi_{12}a_{z_\ell}+\phi_{22}b_{z_\ell}=0,\\
\phi_{12}b_{z_\ell}+\phi_{22}c_{z_\ell}=0,
\end{array}
\end{eqnarray}

If $\phi_{22}\neq 0$ on a non-empty open set, which may happen in all cases covered by Theorems \ref{teo7.2}-\ref{teo7.5}, we obtain from \eqref{31} that
$$
b_{z_\ell} = -\frac{\phi_{12}}{\phi_{22}}a_{z_\ell}, \quad c_{z_\ell} = \left(\frac{\phi_{12}}{\phi_{22}}\right)^2 a_{z_\ell}.
$$
Differentiating the Gauss equation \eqref{gauss} with respect to $z_\ell$ leads to $a_{z_\ell}c+ac_{z_\ell}-2bb_{z_\ell}=0$. Which implies using \eqref{31} that
\begin{eqnarray}\label{33a}
\left[ c + \left(\frac{\phi_{12}}{\phi_{22}}\right)^2 a+2\frac{\phi_{12}}{\phi_{22}} b\right]a_{z_\ell}=0,
\end{eqnarray}

If the expression between brackets in \eqref{33a} does not vanish on a open set, we obtain $a_{z_\ell}=0$ and, thus, by \eqref{31}, $b_{z_\ell}=c_{z_\ell}=0$. Successive differentiation of \eqref{20}, \eqref{21} and \eqref{gauss} with respect to $z_{\ell}, \ldots, z_3$ leads to $a_{z_\ell} = a_{z_{\ell-1}}=\ldots=a_{z_2} = 0$ and, thus, $b_{z_\ell} = b_{z_{\ell-1}}=\ldots=b_{z_2} = 0$ and $c_{z_\ell} = c_{z_{\ell-1}}=\ldots=c_{z_2} = 0$. Therefore, equations \eqref{20} and \eqref{21} give us, respectively,
\begin{eqnarray}\label{20.1}
f_{11}a_t + f_{21}b_t - f_{12}a_x - f_{22}b_x - 2b(f_{11}f_{32}-f_{12}f_{31})+(a-c)(f_{21}f_{32}-f_{22}f_{31})\hspace*{2 cm}\nonumber\\
-\sum_{i=0}^1 (f_{12}a_{z_i}+f_{22}b_{z_i})z_{i+1}+(f_{11}a_{w_0}+f_{21}b_{w_0})w_1 + (f_{11}a_{v_0}+f_{21}b_{v_0})v_1=0,
\end{eqnarray}
and
\begin{eqnarray}\label{21.1}
f_{11}b_t + f_{21}c_t - f_{12}b_x - f_{22}c_x + 2b(f_{21}f_{32}-f_{22}f_{31})+(a-c)(f_{11}f_{32}-f_{12}f_{31})\hspace*{2 cm}\nonumber\\
-\sum_{i=0}^1 (f_{12}b_{z_i}+f_{22}c_{z_i})z_{i+1}+(f_{11}b_{w_0}+f_{21}c_{w_0})w_1 + (f_{11}b_{v_0}+f_{21}c_{v_0})v_1=0.
\end{eqnarray}
Differentiating \eqref{20.1} and \eqref{21.1} with respect to $w_1$, we obtain
\begin{eqnarray}\label{20.2}
f_{11}a_{w_0}+f_{21}b_{w_0} = 0, \quad f_{11}b_{w_0}+f_{21}c_{w_0}=0.
\end{eqnarray}
Likewise,
\begin{eqnarray}\label{21.2}
f_{11}a_{v_0}+f_{21}b_{v_0} = 0, \quad f_{11}b_{v_0}+f_{21}c_{v_0}=0.
\end{eqnarray}
Differentiating the Gauss equation with respect to $w_0$ and $v_0$  leads to $a_{w_0}c+ac_{w_0}-2bb_{w_0}=0$ and $a_{v_0}c+ac_{v_0}-2bb_{v_0}=0$, respectively. Taking into account \eqref{20.2} and \eqref{21.2} in the latter, we obtain
\begin{eqnarray*}
\left[ c + \left(\frac{f_{11}}{f_{21}}\right)^2 a+2\frac{f_{11}}{f_{21}} b\right]a_{w_0}=0,\quad \left[ c + \left(\frac{f_{11}}{f_{21}}\right)^2 a+2\frac{f_{11}}{f_{21}} b\right]a_{v_0}=0,
\end{eqnarray*}
and by \eqref{19.1} we finally have $a_{w_0} = a_{v_0} = 0$ and, thus, by \eqref{20.2} and \eqref{21.2} $b_{w_0} = b_{v_0} = 0$ and $c_{w_0} = c_{v_0} = 0$. Hence, $a$, $b$ and $c$ are universal.

On the other hand, if the expression between brackets in \eqref{33a} vanishes, i.e.,
\begin{eqnarray}\label{5.21a}
c + \left(\frac{\phi_{12}}{\phi_{22}}\right)^2 a+2\frac{\phi_{12}}{\phi_{22}} b=0,
\end{eqnarray}
then it follows from \eqref{5.21a} and \eqref{gauss} that
\begin{eqnarray}\label{35}
 b = \pm 1 - \frac{\phi_{12}}{\phi_{22}}a,\quad  c = \left( \frac{\phi_{12}}{\phi_{22}} \right)^2 a \mp 2 \frac{\phi_{12}}{\phi_{22}}.
\end{eqnarray}
Therefore, 
\begin{eqnarray*}
f_{11}D_t a+f_{21}D_t b &=& \frac{\Delta_{12}}{\phi_{22}}D_ta - af_{21}\left( \frac{\phi_{12}}{\phi_{22}} \right)_{,t},\\
f_{12}D_x a + f_{22}D_x b &=& -\lambda z_0 \frac{\Delta_{12}}{\phi_{22}}D_xa - af_{22}\left( \frac{\phi_{12}}{\phi_{22}} \right)_{,x},\\
f_{11}D_t b+f_{21}D_t c &=& \left( 2f_{21}\frac{\phi_{12}}{\phi_{22}}a-f_{11}a\mp 2f_{21}\right)\left( \frac{\phi_{12}}{\phi_{22}} \right)_{,t}-\frac{\phi_{12}}{\phi_{22}^2}\Delta_{12}D_ta ,\\
f_{12}D_x b + f_{22}D_x c &=& \left[\left( \lambda z_0\frac{\Delta_{12}}{\phi_{22}}+f_{22}\frac{\phi_{12}}{\phi_{22}}\right)a \mp 2f_{22}\right]\left( \frac{\phi_{12}}{\phi_{22}} \right)_{,x}+\lambda z_0 \frac{\phi_{12}}{\phi_{22}^2}\Delta_{12}D_xa,
\end{eqnarray*}
where $\Delta_{12} = f_{11}\phi_{22}-f_{21}\phi_{12}$.

Therefore, equation \eqref{3.8} becomes
{\small\begin{eqnarray}\label{42}
\frac{\Delta_{12}}{\phi_{22}}(D_ta+\lambda z_0D_xa) - af_{21}\left( \frac{\phi_{12}}{\phi_{22}} \right)_{,t}  + af_{22}\left( \frac{\phi_{12}}{\phi_{22}} \right)_{,x}-2b\Delta_{13}+(a-c)\Delta_{23} = 0,
\end{eqnarray}}
and \eqref{3.9} becomes
{\small\begin{eqnarray}\label{43}
-\frac{\phi_{12}}{\phi_{22}^2}\Delta_{12}(D_ta+\lambda z_0D_xa)+\left( 2f_{21}\frac{\phi_{12}}{\phi_{22}}a-f_{11}a\mp 2f_{21}\right)\left( \frac{\phi_{12}}{\phi_{22}} \right)_{,t}\nonumber\hspace*{3 cm}\\
-\left[\left( \lambda z_0\frac{\Delta_{12}}{\phi_{22}}+f_{22}\frac{\phi_{12}}{\phi_{22}}\right)a \mp 2f_{22}\right]\left( \frac{\phi_{12}}{\phi_{22}} \right)_{,x}+(a-c)\Delta_{13}+2b\Delta_{23} = 0,
\end{eqnarray}}
where
$$
\left( \frac{\phi_{12}}{\phi_{22}} \right)_{,t}=\left[\left( \frac{\phi_{12}}{\phi_{22}} \right)_{,z_0}w_1 + \left( \frac{\phi_{12}}{\phi_{22}} \right)_{,z_1}v_1 \right].
$$
From Lemma \ref{lemma3.1} we have $\phi_{12}\neq 0$, since $c\neq 0$. Hence, adding \eqref{42} multiplied by $\phi_{12}/\phi_{22}$ to \eqref{43} we get
{\small\begin{eqnarray}\label{44}
\left( -\frac{\Delta_{12}}{\phi_{22}}a\mp 2f_{21}\right)\left( \frac{\phi_{12}}{\phi_{22}} \right)_{,t}+\left( -\lambda z_0\frac{\Delta_{12}}{\phi_{22}}a \pm 2f_{22}\right) \left( \frac{\phi_{12}}{\phi_{22}} \right)_{,x}+\left(a-2b\frac{\phi_{12}}{\phi_{22}}-c\right)\Delta_{13}\nonumber\hspace*{2 cm}\\
+\left[\frac{\phi_{12}}{\phi_{22}}(a-c)+2b\right]\Delta_{23} = 0.
\end{eqnarray}}
Differentiating \eqref{44} with respect to $v_1$ and $w_1$, we obtain, respectively,
\begin{eqnarray}\label{4.17.1}
P\left( \frac{\phi_{12}}{\phi_{22}} \right)_{,z_1} = 0,\quad P\left( \frac{\phi_{12}}{\phi_{22}} \right)_{,z_0} = 0,\quad \textnormal{ where }\quad P:= -\frac{\Delta_{12}}{\phi_{22}}a\mp 2f_{21}.
\end{eqnarray}
If $P\neq 0$ on a non-empty open set, we have from \eqref{4.17.1} that $\phi_{22} - A\phi_{12} = 0$, where $A$ is a nonzero constant. But, $\ell = \phi_{22} - A\phi_{12} = 0$ restricts our analysis to the case where $f_{ij}$ are given by \eqref{th2.2a} with $A=\mu\neq 0$ or \eqref{th2.3a} with $A=\mu\neq 0$.

If $f_{ij}$ are given by \eqref{th2.2a} with $\mu\neq 0$ we obtain $\Delta_{13} = f_{11}f_{32}-f_{31}f_{12} = \mp m\mu\psi$ and $\Delta_{23} = f_{21}f_{32}-f_{31}f_{22} = \pm m\psi$, which imply by \eqref{44} that $a = \pm 2\mu / (1+\mu^2)$ constant. Therefore, \eqref{42} and \eqref{43} reduce to
\begin{equation*}
\left(\begin{array}{ccc}
-2b & a-c \\  a-c & 2b
\end{array}\right)
\left(\begin{array}{cc}
\Delta_{13} \\ \Delta_{23}
\end{array}\right)
=
\left(\begin{array}{cc}
0 \\ 0
\end{array}\right).
\end{equation*}
It follows from \eqref{3.10a} that $b=0$ and $a=c$, which contradicts the Gauss equation \eqref{gauss}.

If $f_{ij}$ are given by \eqref{th2.3a} with $\mu\neq 0$ we obtain $\Delta_{13} = f_{11}f_{32}-f_{31}f_{12} = \lambda (m_1\mu -\eta)z_1$ and $\Delta_{23} = f_{21}f_{32}-f_{31}f_{22} = -\lambda m_1z_1$, which imply by \eqref{44} that 
\begin{eqnarray}\label{4.18.1}
\left[ a \frac{1+\mu^2}{\mu}\mp 2\right]m_1 = 0.
\end{eqnarray}
Therefore, if $m_1\neq 0$ in \eqref{4.18.1} then $a$ is a constant and \eqref{42} and \eqref{43} give us a contradiction like before. If $m_1 = 0$ in \eqref{4.18.1} then replacing $\Delta_{13} = -\lambda \eta z_1 (\neq 0)$ and $\Delta_{23}=0$ into \eqref{42} and \eqref{43}, we get
\begin{eqnarray*}
D_t a+\lambda z_0D_xa + 2\lambda(\pm \mu - a)z_1 = 0,&&\\
D_t a+\lambda z_0D_xa + \lambda (a\mu^2-a\pm 2\mu)z_1 = 0,&&
\end{eqnarray*}
which imply $a=0$ and, therefore, a contradiction by Lemma \ref{lemma3.1}.

On the other hand, from \eqref{4.17.1} if $P = 0$, on a open set, then $a = \mp 2\phi_{22}f_{21}/\Delta_{12}$ and, thus, $a$ is a function depending only on $z_0$, $z_1$ and $z_2$. However, using such $a$ and \eqref{35} we get
$$
c + \left(\frac{f_{11}}{f_{21}}\right)^2 a+2\frac{f_{11}}{f_{21}} b = 0,
$$
which contradicts the hypothesis \eqref{19.1}.

If $\phi_{22} \equiv 0$ (which is the case of \eqref{th2.2a} with $\mu=0$ or \eqref{th2.3a} with $\mu=0$ or \eqref{th2.4a} with $\mu\neq 0$ and $\psi = -(m_2+\lambda m_1z_0)\sqrt{1+\mu^2}/\mu$) then,  since $\Delta_{12} = [(\phi_{22}-\mu_2\phi_{12})f_{11}-\eta_2\phi_{12}]dx\wedge dt = -\phi_{12}f_{21}dx\wedge dt\neq 0$, we have $\phi_{12}\neq 0$ on an open set. Moreover, it follows from \eqref{31} that $a_{z_l} = b_{z_l} = 0$. 

Differentiating the Gauss equation with respect to $z_\ell$ and using the Lemma \eqref{lemma3.1}, we obtain $c_{z_\ell}=0$. Successive differentiation of \eqref{20}, \eqref{21} and \eqref{gauss} with respect to $z_l, \ldots, z_3$ leads to $a_{z_i}=b_{z_i}=c_{z_i} = 0$ for $i = 2,3,...,l-1$. Hence, \eqref{20} and \eqref{21} are equivalent to
\begin{eqnarray}\label{20a}
f_{11}a_t + f_{21}b_t - f_{12}a_x - f_{22}b_x - 2b(f_{11}f_{32}-f_{12}f_{31})+(a-c)(f_{21}f_{32}-f_{22}f_{31})\hspace*{2 cm}\nonumber\\
-\sum_{i=0}^1 (f_{12}a_{z_i}+f_{22}b_{z_i})z_{i+1}+(f_{11}a_{w_0}+f_{21}b_{w_0})w_1 + (f_{11}a_{v_0}+f_{21}b_{v_0})v_1=0,
\end{eqnarray}
and
\begin{eqnarray}\label{21a}
f_{11}b_t + f_{21}c_t - f_{12}b_x - f_{22}c_x + 2b(f_{21}f_{32}-f_{22}f_{31})+(a-c)(f_{11}f_{32}-f_{12}f_{31})\hspace*{2 cm}\nonumber\\
-\sum_{i=0}^1 (f_{12}b_{z_i}+f_{22}c_{z_i})z_{i+1}+(f_{11}b_{w_0}+f_{21}c_{w_0})w_1 + (f_{11}b_{v_0}+f_{21}c_{v_0})v_1=0.
\end{eqnarray}
Differentiating \eqref{20a} and \eqref{21a} with respecto to $v_1$ and $w_1$, we get
\begin{eqnarray}\label{20.2a}
\begin{array}{rr}
f_{11}a_{w_0}+f_{21}b_{w_0} = 0,\\ \quad f_{11}b_{w_0}+f_{21}c_{w_0}=0,
\end{array}
\begin{array}{rr}
f_{11}a_{v_0}+f_{21}b_{v_0} = 0,\\ \quad f_{11}b_{v_0}+f_{21}c_{v_0}=0.
\end{array}
\end{eqnarray}
Differentiating the Gauss equation \eqref{gauss} with respect to $w_0$ and $v_0$  leads to $a_{w_0}c+ac_{w_0}-2bb_{w_0}=0$ and $a_{v_0}c+ac_{v_0}-2bb_{v_0}=0$, respectively. Taking into account \eqref{20.2a} in the latter, by \eqref{19.1} we obtain $a_{w_0} = a_{v_0} = 0$ and, thus, by \eqref{20.2a}, $b_{w_0} = b_{v_0} = 0$ and $c_{w_0} = c_{v_0} = 0$. Hence, $a$, $b$ and $c$ are universal and, thus, we conclude the proof of $(i)$.

Suppose $(ii)$, i.e, $m < n$. Therefore, since $n\geq m+1$, differentiating \eqref{10}, \eqref{11} and \eqref{gauss} with respect to $v_{n+1}$ leads to $a_{v_{n}}=b_{v_{n}}=c_{v_{n}}=0$. Successive differentiation with respect to $v_n, v_{n-1}, \ldots, v_{(m+1)+1}$ leads to $a_{v_{n-1}}=\ldots = a_{v_{m+1}} = 0$, $b_{v_{n-1}}=\ldots = b_{v_{m+1}} = 0$ and $c_{v_{n-1}}=\ldots = c_{v_{m+1}} = 0$. Hence, $a$, $b$ and $c$ are functions of $x, t, z_0, z_1, \ldots, z_l, w_1, \ldots, w_m$, $v_1, \ldots, v_m$. Proceeding as in $(i)$, we conclude that $a$, $b$ and $c$ are functions of $x$ and $t$ only, and therefore universal. This concludes $(ii)$.

Finally $(iii)$, i.e, $m > n$. Therefore, since $m\geq n+1$, differentiating \eqref{10}, \eqref{11} and \eqref{gauss} with respect to $w_{m+1}$ leads to $a_{w_{m}}=b_{w_{m}}=c_{w_{m}}=0$. Successive differentiation with respect to $w_m, w_{m-1}, \ldots, w_{(n+1)+1}$ leads to $a_{w_{m-1}}=\ldots = a_{w_{n+1}} = 0$, $b_{w_{m-1}}=\ldots = b_{w_{n+1}} = 0$ and $c_{w_{m-1}}=\ldots = c_{w_{n+1}} = 0$. Hence, $a$, $b$ and $c$ are functions of $x, t, z_0, z_1, \ldots, z_l, w_1, \ldots, w_n, v_1, \ldots, v_n$. Proceeding as in $(i)$, we conclude that $a$, $b$ and $c$ are functions of $x$ and $t$ only, and therefore universal. This concludes $(iii)$.

Therefore, $a$, $b$ and $c$ are universal, i.e., $a$, $b$ and $c$ depend only on $x$ and $t$. This concludes the proof of Lemma \ref{lemma3.2}.

\cqd


\begin{lemma}\label{lemma3.4}
Consider an equation of type \eqref{T} describing pseudospherical surfaces, under the condition \eqref{1}, given by the Theorems $\ref{teo7.2}$-$\ref{teo7.5}$. Assume there is a local isometric immersion of the pseudospherical surface, determined by a solution $u(x,t)$ of \eqref{T}, for which the coefficients $a$, $b$ and $c$ of the second fundamental form depend on $x, t, z_0, \ldots, z_l, w_1, \ldots, w_m, v_1, \ldots, v_n$, where $1\leq l  <\infty$, $1\leq m < \infty$ and $1\leq n < \infty$ are finite, but otherwise arbitrary. If $f_{21} = 0$, on a non-empty open set, then $a$, $b$ and $c$ are functions of $x$ and $t$ only, and therefore universal.
\end{lemma}

\noindent\textbf{Proof}. First, observe that $f_{21} = 0$, on a non-empty open set, can only happen in \eqref{th2.4a} with $\mu=m_1 = 0$ or \eqref{th2.5b} with $\mu = \eta = 0$. Furthermore, in both cases $\phi_{22}\neq 0$ on that open set. Our analysis consists in three cases, namely,
$$
(i)\hspace{0.2 cm} m=n, \qquad (ii)\hspace{0.2 cm} m < n, \qquad (iii)\hspace{0.2 cm} n < m.
$$

Let consider, firstly, the case $m=n$. Suppose $l =1$. Successive differentiation of \eqref{10}, \eqref{11} and \eqref{gauss} with respect to $v_{n+1}, \ldots, v_1$ and $w_{n+1}, \ldots, w_1$, since $f_{11}\neq 0$, lead to $a_{v_k}=b_{v_k}=c_{v_k}=0$ and $a_{w_k}=b_{w_k}=c_{w_k}=0$ for $k=0,1, \ldots, n$. Therefore, $a$, $b$ and $c$ are universal. 

Now let us consider $l \geq 2$. Taking successive differentiation of \eqref{10}, \eqref{11} and \eqref{gauss} with respect to $v_{n+1}, \ldots, v_2$ and $w_{n+1}, \ldots, w_2$, since $f_{11}\neq 0$, leads to $a_{w_k}=b_{w_k}=c_{w_k}=0$ and $a_{v_k}=b_{v_k}=c_{v_k}=0$ for $k=1,2, \ldots, n$. Thus, we have that $a$, $b$ and $c$ do not depend on $w_k$ and neither $v_k$ for $k=1,2, \ldots, n$. Hence, $a$, $b$ and $c$ are functions of $x, t, z_0=w_0, z_1=v_0, \ldots, z_l$. Furthermore, the equations \eqref{10} and \eqref{11} are equivalent to
\begin{eqnarray}\label{a20}
f_{11}a_t  - f_{12}a_x - f_{22}b_x - 2b(f_{11}f_{32}-f_{12}f_{31})-(a-c)f_{22}f_{31}-\sum_{i=0}^l (f_{12}a_{z_i}+f_{22}b_{z_i})z_{i+1}\nonumber\hspace*{1 cm}\\
+\sum_{i=2}^l f_{11}a_{z_i}\partial_x^{i-2}(z_{0,t}-F)+f_{11}a_{w_0}w_1 + f_{11}a_{v_0}v_1=0,
\end{eqnarray}
and
\begin{eqnarray}\label{a21}
f_{11}b_t  - f_{12}b_x - f_{22}c_x - 2bf_{22}f_{31}+(a-c)(f_{11}f_{32}-f_{12}f_{31})-\sum_{i=0}^l (f_{12}b_{z_i}+f_{22}c_{z_i})z_{i+1}\nonumber\hspace*{1 cm}\\
+\sum_{i=2}^l f_{11}b_{z_i}\partial_x^{i-2}(z_{0,t}-F)+f_{11}b_{w_0}w_1 + f_{11}b_{v_0}v_1=0.
\end{eqnarray}

Differentiating \eqref{a20} and \eqref{a21} with respect to $z_{l+1}$, we obtain, respectively,
\begin{eqnarray*}
(f_{12}+\lambda z_0f_{11})a_{z_l}+f_{22}b_{z_l}=0,\quad (f_{12}+\lambda z_0f_{11})b_{z_l}+f_{22}c_{z_l}=0,
\end{eqnarray*}
and, using \eqref{(7.2)}, with $f_{21}=0$, we have
\begin{eqnarray}
\begin{array}{rr}\label{a31}
\phi_{12}a_{z_l}+\phi_{22}b_{z_l}=0,\quad \phi_{12}b_{z_l}+\phi_{22}c_{z_l}=0,
\end{array}
\end{eqnarray}

Since $\phi_{22}\neq 0$, it follows from \eqref{a31} that
\begin{eqnarray}\label{(90)}
b_{z_l} = -\frac{\phi_{12}}{\phi_{22}}a_{z_l}, \quad c_{z_l} = \left(\frac{\phi_{12}}{\phi_{22}}\right)^2 a_{z_l}.
\end{eqnarray}
Differentiating the Gauss equation \eqref{gauss} with respect to $z_l$ leads to $a_{z_l}c+ac_{z_l}-2bb_{z_l}=0$, which gives, using the \eqref{(90)}
\begin{eqnarray}\label{a33}
\left[ c + \left(\frac{\phi_{12}}{\phi_{22}}\right)^2 a+2\frac{\phi_{12}}{\phi_{22}} b\right]a_{z_l}=0,
\end{eqnarray}

If the expression between brackets in \eqref{a33} does not vanish on a open set, we obtain $a_{z_l}=0$ and, thus, by \eqref{a31}, $b_{z_l}=c_{z_l}=0$. Successive differentiation of \eqref{a20}, \eqref{a21} and \eqref{gauss} with respect to $z_{l}, \ldots, z_3$ leads to $a_{z_l} = a_{z_{l-1}}=\ldots=a_{z_2} = 0$ and, thus, $b_{z_l} = b_{z_{l-1}}=\ldots=b_{z_2} = 0$ and $c_{z_l} = c_{z_{l-1}}=\ldots=c_{z_2} = 0$. Therefore, the equation \eqref{a20} and \eqref{a21} give us, respectively,
\begin{eqnarray}\label{a20.1}
f_{11}a_t  - f_{12}a_x - f_{22}b_x - 2b(f_{11}f_{32}-f_{12}f_{31})-(a-c)f_{22}f_{31}-\sum_{i=0}^2 (f_{12}a_{z_i}+f_{22}b_{z_i})z_{i+1}\nonumber\hspace*{1 cm}\\
+f_{11}a_{w_0}w_1 + f_{11}a_{v_0}v_1=0,
\end{eqnarray}
and
\begin{eqnarray}\label{a21.1}
f_{11}b_t  - f_{12}b_x - f_{22}c_x - 2bf_{22}f_{31}+(a-c)(f_{11}f_{32}-f_{12}f_{31})-\sum_{i=0}^2 (f_{12}b_{z_i}+f_{22}c_{z_i})z_{i+1}\nonumber\hspace*{1 cm}\\
+f_{11}b_{w_0}w_1 + f_{11}b_{v_0}v_1=0.
\end{eqnarray}
Differentiating \eqref{a20.1} and \eqref{a21.1} with respect to $v_1$ and $w_1$ leads to $f_{11}a_{v_0} = f_{11}b_{v_0}=0$ and $f_{11}a_{w_0}= f_{11}b_{w_0}=0$, i.e, $a_{v_0} =b_{v_0} =0$ and $a_{w_0}=b_{w_0}=0$. Differentiating the Gauss equation \eqref{gauss} with respect to $w_0$ and $v_0$  gives $a_{w_0}c+ac_{w_0}-2bb_{w_0}=0$ and $a_{v_0}c+ac_{v_0}-2bb_{v_0}=0$, respectively. Since $a\neq 0$ we obtain $c_{w_0} = c_{v_0} = 0$. Hence, $a$, $b$ and $c$ are universal.

On the other hand, if the expression in brackets in \eqref{a33} vanishes, i.e.,
\begin{eqnarray}\label{a5.14.1}
c + \left(\frac{\phi_{12}}{\phi_{22}}\right)^2 a+2\frac{\phi_{12}}{\phi_{22}} b=0,
\end{eqnarray}
then it follows from the Gauss equation \eqref{gauss} that
\begin{eqnarray}
&& b = \pm 1 - \frac{\phi_{12}}{\phi_{22}}a,\label{a5.15}\\
&& c = \left( \frac{\phi_{12}}{\phi_{22}} \right)^2 a \mp 2 \frac{\phi_{12}}{\phi_{22}}.\label{a5.16}
\end{eqnarray}

Therefore, 
\begin{eqnarray*}
f_{12}D_x a + f_{22}D_x b &=& -\lambda z_0 f_{11}D_xa - af_{22}\left( \frac{\phi_{12}}{\phi_{22}} \right)_{,x},\\
f_{11}D_t b &=& -f_{11}a\left( \frac{\phi_{12}}{\phi_{22}} \right)_{,t}-\frac{\phi_{12}}{\phi_{22}}f_{11}D_ta ,\\
f_{12}D_x b + f_{22}D_x c &=& \left[\left( \lambda z_0f_{11}+\phi_{12}\right)a \mp 2f_{22}\right]\left( \frac{\phi_{12}}{\phi_{22}} \right)_{,x}+\lambda z_0 \frac{\phi_{12}}{\phi_{22}}f_{11}D_xa,
\end{eqnarray*}
where $\Delta_{12} = f_{11}\phi_{22}-f_{21}\phi_{12}$.

Therefore, equation \eqref{3.8} becomes
{\small\begin{eqnarray}\label{a42}
f_{11}(D_ta+\lambda z_0D_xa)  + af_{22}\left( \frac{\phi_{12}}{\phi_{22}} \right)_{,x}-2b\Delta_{13}+(a-c)\Delta_{23} = 0,
\end{eqnarray}}
and \eqref{3.9} becomes
{\small\begin{eqnarray}\label{a43}
-\frac{\phi_{12}}{\phi_{22}}f_{11}(D_ta+\lambda z_0D_xa)-f_{11}a\left( \frac{\phi_{12}}{\phi_{22}} \right)_{,t}-\left[\left( \lambda z_0f_{11}+\phi_{12}\right)a \mp 2f_{22}\right]\left( \frac{\phi_{12}}{\phi_{22}} \right)_{,x}+(a-c)\Delta_{13}+2b\Delta_{23} = 0,
\end{eqnarray}}
where
$$
\left( \frac{\phi_{12}}{\phi_{22}} \right)_{,t}=\left[\left( \frac{\phi_{12}}{\phi_{22}} \right)_{,z_0}w_1 + \left( \frac{\phi_{12}}{\phi_{22}} \right)_{,z_1}v_1 \right].
$$
From Lemma \ref{lemma3.1}, since $c\neq 0$ we have $\phi_{12}\neq 0$. Hence, adding \eqref{a42} multiplied by $\phi_{12}/\phi_{22}$ with \eqref{a43} we get
{\small\begin{eqnarray}\label{a44}
 -\frac{\phi_{12}}{\phi_{22}}f_{11}a\left( \frac{\phi_{12}}{\phi_{22}} \right)_{,t}+\left( -\lambda z_0f_{11}a \pm 2f_{22}\right) \left( \frac{\phi_{12}}{\phi_{22}} \right)_{,x}+\left(a-2b\frac{\phi_{12}}{\phi_{22}}-c\right)\Delta_{13}+\left[\frac{\phi_{12}}{\phi_{22}}(a-c)+2b\right]\Delta_{23} = 0.
\end{eqnarray}}
Differentiating \eqref{a44} with respect to $v_1$ and $w_1$, we obtain, respectively,
\begin{eqnarray*}
 -\frac{\phi_{12}}{\phi_{22}}f_{11}a\left( \frac{\phi_{12}}{\phi_{22}} \right)_{,z_1} = 0,\quad  -\frac{\phi_{12}}{\phi_{22}}f_{11}a\left( \frac{\phi_{12}}{\phi_{22}} \right)_{,z_0} = 0,
\end{eqnarray*}
which imply that $\phi_{22} - A\phi_{12} = 0$, where $A$ is a nonzero constant. Otherwise, we would have $\phi_{22}=0$. But, $l = \phi_{22} - A\phi_{12} = 0$ does not happen in \eqref{th2.4a} or \eqref{th2.5b}. This concludes $(i)$.

Suppose $(ii)$, i.e, the case $m < n$. Therefore, since $n\geq m+1$, differentiating \eqref{10}, \eqref{11} and \eqref{gauss} with respect to $v_{n+1}$ leads to $a_{v_{n}}=b_{v_{n}}=c_{v_{n}}=0$. Successive differentiation with respect to $v_n, v_{n-1}, \ldots, v_{(m+1)+1}$ leads to $a_{v_{n-1}}=\ldots = a_{v_{m+1}} = 0$, $b_{v_{n-1}}=\ldots = b_{v_{m+1}} = 0$ and $c_{v_{n-1}}=\ldots = c_{v_{m+1}} = 0$. Hence, $a$, $b$ and $c$ are functions of $x, t, z_0, z_1, \ldots, z_l, w_1, \ldots, w_m$, $v_1, \ldots, v_m$. Proceeding as in $(i)$, we conclude that $a$, $b$ and $c$ are functions of $x$ and $t$ only, and thus universal. This concludes $(ii)$.

To conclude the proof of Lemma \ref{lemma3.4}, consider $(iii)$, i.e, the case $m > n$. Therefore, since $m\geq n+1$, differentiating \eqref{10}, \eqref{11} and \eqref{gauss} with respect to $w_{m+1}$ leads to $a_{w_{m}}=b_{w_{m}}=c_{w_{m}}=0$. Successive differentiation with respect to $w_m, w_{m-1}$, $\ldots$, $w_{(n+1)+1}$ leads to $a_{w_{m-1}}=\ldots = a_{w_{n+1}} = 0$, $b_{w_{m-1}}=\ldots = b_{w_{n+1}} = 0$ and $c_{w_{m-1}}=\ldots = c_{w_{n+1}} = 0$. Hence, $a$, $b$ and $c$ are functions of $x, t$, $z_0, z_1, \ldots, z_l, w_1, \ldots, w_n, v_1, \ldots, v_n$. Again, proceeding as in $(i)$, we conclude that $a$, $b$ and $c$ are functions of $x$ and $t$ only, and thus universal. This concludes $(iii)$.

Hence, $a$, $b$ and $c$ are universal, i.e., $a$, $b$ and $c$ depend only on $x$ and $t$. This concludes the proof of Lemma \ref{lemma3.4}.

\cqd

\subsection{Universal expressions for the second fundamental forms}\label{sec5}

In the previous section we have shown that if there exist coefficients $a$, $b$, $c$ (depending on a jet of finite order of $u$) of the second fundamental form of a local isometric immersion of a pseudospherical surface, so that the system of equations \eqref{gauss}, \eqref{3.8} and \eqref{3.9} is satisfied, then $a$, $b$ and $c$ are functions depending only on $x$ and $t$, and thus universal. Now we are going to determine such coefficients for the equations \eqref{eqT} and associated $f_{ij's}$ given by Theorems \ref{teo7.2}-\ref{teo7.5}.

\begin{proposition}\label{prop1}
Consider an equation of type \eqref{T} describing pseudospherical surfaces, under the condition \eqref{1}, given by the Theorem $\ref{teo7.2}$. There exists a local isometric immersion in $\mathbb{R}^3$ of a pseudospherical surface, defined by a solution $u$, for which the coefficients $a$, $b$ and $c$ of the second fundamental form depend on $x, t, z_0, \ldots, z_l, w_1, \ldots, w_m$, $v_1, \ldots, v_n$, where $1\leq l  <\infty$, $1\leq m < \infty$ and $1\leq n < \infty$ are finite, but otherwise arbitrary, if, and only if, 

\vspace*{0.4 cm}
\noindent
$(i)$ $\mu =0$ and $a$, $b$ and $c$ depend only on $x$ and are given by
\begin{eqnarray}\label{7.1b}
a = \pm \sqrt{L(x)}, \quad b=-\beta e^{\pm 2\eta x}, \quad c = a \mp \frac{a_x}{\eta},
\end{eqnarray}
where $L(x) = \sigma e^{\pm 2\eta x}-\beta^2 e^{\pm 4\eta x}-1$, with $\eta$, $\sigma$, $\beta$ $\in$ $\mathbb{R}$, $\eta\neq 0$, $\sigma >0$ and $\sigma^2 >4\beta^2$. The coefficients $a$, $b$, $c$ are defined on a strip of $\mathbb{R}$ where
\begin{eqnarray}\label{p7.2}
log\sqrt{\frac{\sigma -\sqrt{\sigma^2 -4\beta^2}}{2\beta^2}}< \pm \eta x <log\sqrt{\frac{\sigma +\sqrt{\sigma^2 -4\beta^2}}{2\beta^2}}.
\end{eqnarray}
Moreover, the constants $\beta$ and $\sigma$ have to be chosen so that the strip intersects the domain of the solution of \eqref{th2.2}.

\vspace*{0.3 cm}
\noindent
or
\vspace*{0.3 cm}

\noindent
$(ii)$ $\mu \neq 0$ and $a$, $b$ and $c$ depend only on $x$ and are given by
\begin{eqnarray}\label{101.a}
 a &=& \frac{1}{2\mu}[\pm \mu\sqrt{\Delta}-(\mu^2-1)b+\beta e^{\pm 2\eta x}],\nonumber\\
 c &=& \frac{1}{2\mu}[\pm \mu\sqrt{\Delta}+(\mu^2-1)b-\beta e^{\pm 2\eta x}],\\
 \Delta &=& \frac{[(\mu^2-1)b-\beta e^{\pm 2\eta x}]^2 -4\mu^2(1-b^2)}{\mu^2}> 0\nonumber
\end{eqnarray}
where $b$ satisfies the ordinary differential equation
\begin{eqnarray}\label{102}
[\mu(1+\mu^2)\sqrt{\Delta}\pm (\mu^2+1)^2b\mp (\mu^2-1)\beta e^{\pm 2\eta x}]b'\hspace*{3cm}\nonumber\\
+ 2\eta\left\lbrace [\mp\mu (1+\mu^2)\sqrt{\Delta}- \beta(\mu^2-1)e^{\pm 2\eta x}]b+ \beta^2e^{\pm 4\eta x}\right\rbrace=0
\end{eqnarray}
\end{proposition}

\noindent
\textbf{Proof}. Since $\eta=m\neq 0$ we only have $f_{21}\neq 0$, on a open set. From Lemma \ref{lemma3.3} the equations \eqref{gauss}, \eqref{3.8} and \eqref{3.9} form an inconsistent system. From Lemma \ref{lemma3.2}, the coefficients of the second fundamental form of such local isometric immersion are universal, and hence \eqref{3.8} and \eqref{3.9} become
\begin{eqnarray}
f_{11} a_t+f_{21} b_t -f_{12} a_x-f_{22} b_x-2b\Delta_{13}+(a-c)\Delta_{23} = 0,&&\label{7.1}\\
f_{11} b_t+f_{21} c_t -f_{12} b_x-f_{22} c_x+(a-c)\Delta_{13}+2b\Delta_{23} = 0,&&\label{7.2}
\end{eqnarray}
where $\Delta_{13} = \mp \eta\mu\psi$ and $\Delta_{23}=\pm \eta\psi (\neq 0)$. Hence, since $f_{ij}$ are given by \eqref{th2.2a}, differentiating \eqref{7.1} and \eqref{7.2} with respect to $z_2$ we obtain
\begin{eqnarray}
&&a_t = \mu^2 c_t,\label{7.3}\\
&&b_t = -\mu c_t.\label{7.4}
\end{eqnarray}
Replacing \eqref{7.3} and \eqref{7.4} back into \eqref{7.1} and \eqref{7.2} we get
\begin{eqnarray}
-\mu \eta\sqrt{1+\mu^2} c_t + [-a_x-\mu b_x \pm \eta(2\mu b + a-c)]\psi = 0,&&\label{7.5}\\
\eta\sqrt{1+\mu^2}c_t + [-b_x-\mu c_x\pm \eta(2b - \mu a+\mu c)]\psi =0.&&\label{7.6}
\end{eqnarray}
Isolating $\eta\sqrt{1+\mu^2}c_t$ in \eqref{7.6} and replacing it into \eqref{7.5}, we obtain 
\begin{eqnarray}\label{7.7}
\mu [-b_x-\mu c_x\pm \eta(2b - \mu a+\mu c)] + [-a_x-\mu b_x \pm \eta(2\mu b + a-c)] = 0.
\end{eqnarray}
Differentiating \eqref{7.7} with respect to $t$ and using \eqref{7.3} and \eqref{7.4}, we get $\mp \eta (1+\mu^2)^2c_t=0$ and, thus, $c_t = 0$. By \eqref{7.3} and \eqref{7.4}, $a_t = b_t = 0$. Hence, $a$, $b$ and $c$ are functions depending only on $x$. Therefore, it follows from \eqref{7.5} and \eqref{7.6} that
\begin{eqnarray}
-a_x-\mu b_x \pm \eta(2\mu b + a-c) = 0,&&\label{7.8}\\
-b_x-\mu c_x\pm \eta(2b - \mu a+\mu c) =0,&&\label{7.9}
\end{eqnarray}
where \eqref{7.7} is now identically satisfied. From \eqref{7.8} we have $c$ in terms of $a$, $b$, $a_x$ and $b_x$, which replaced into \eqref{7.9} leads to
\begin{eqnarray}\label{7.12b}
\mu a_x = \pm \eta(1+\mu^2)b - \mu^2 b_x \pm \beta \eta e^{\pm 2\eta x},
\end{eqnarray}
where $\beta$ is a constant.

If $\mu = 0$, then from \eqref{7.12b} and \eqref{7.8}, we have
\begin{eqnarray}\label{7.13b}
b = -\beta e^{\pm 2\eta x}\quad \textnormal{e} \quad c = a \mp \frac{a_x}{\eta}.
\end{eqnarray}
Substituting \eqref{7.13b} in the Gauss equation \eqref{gauss}, we obtain $a = \pm \sqrt{L(x)}$ where $L(x) = \sigma e^{\pm 2\eta x}-\beta^2 e^{\pm 4\eta x}-1$, with $\sigma$, $\beta$ $\in$ $\mathbb{R}$, $\sigma >0$ and $\sigma^2 >4\beta^2$. This $a$ together with \eqref{7.13b} give us \eqref{7.1b}, where $a$ is defined on the strip described by \eqref{p7.2}.

If $\mu \neq 0$ then \eqref{7.12b} gives us $a_x$, which replaced into \eqref{7.8} implies that
\begin{eqnarray}\label{112}
c = a + \phi(x), \quad \phi(x) = \frac{\mu^2-1}{\mu}b(x) - \frac{\beta}{\mu}e^{\pm 2\eta x}.
\end{eqnarray}
Substituting the latter into the Gauss equation we obtain $a^2+a\phi(x)-b^2 =-1$, which resolved as a second degree equation in terms of $a$ leads to
$$
a = \frac{-\phi(x)\pm \sqrt{\Delta}}{2}, \quad \Delta = \phi(x)^2-4[1-b(x)^2]> 0.
$$
Hence, using \eqref{112} we also have $c$ in terms of $b=b(x)$ as in \eqref{101.a}, which replaced into \eqref{7.12b} gives us
\begin{eqnarray}\label{eq}
[(1+\mu^2)\sqrt{\Delta}\pm (\mu^2-1)\phi \pm 4\mu b]b' \mp 2\eta(1+\mu^2)b\sqrt{\Delta}-2\eta\beta e^{\pm 2\eta x}\phi =0.
\end{eqnarray}
Observe that, if the coefficient of $b'$ in \eqref{eq} vanishes, we have
\begin{eqnarray*}
(1+\mu^2)\sqrt{\Delta}\pm (\mu^2-1)\phi \pm 4\mu b=0,&&\\
\mp 2\eta(1+\mu^2)b\sqrt{\Delta}-2\eta\beta e^{\pm 2\eta x}\phi =0.&&
\end{eqnarray*}
In the latter two equations, replacing $(1+\mu^2)\sqrt{\Delta}$ of the first into the second implies that
\begin{eqnarray*}
0 &=& \mp 2\eta b[\mp (\mu^2-1)\phi \mp 4\mu b]-2\eta\beta e^{\pm 2\eta x}\phi\\
  &=& 2\eta\mu (\phi^2 +4b^2),
\end{eqnarray*}
and then $\phi^2 +4b^2=0$, since $\eta\mu\neq 0$. However, $\phi^2 +4b^2=0$ if, and only if, $\phi=b=0$ which implies by \eqref{7.12b} that $a=c$. But, $a=c$ and $b=0$ contradict the Gauss equation \eqref{gauss}. Therefore, the coefficient of $b'$ in the equation \eqref{eq} does not vanish in a non-empty open set. That means we can write $b' = g(x,b)$, where $g$ is a differentiable function defined, from \eqref{eq}, by
$$
g(x,b) = \frac{\pm 2\eta(1+\mu^2)b\sqrt{\Delta}+2\eta\beta e^{\pm 2\eta x}\phi}{(1+\mu^2)\sqrt{\Delta}\pm (\mu^2-1)\phi \pm 4\mu b}.
$$

Let $x_0$ be an arbitrarily fixed point and consider the following Initial Value Problem (IVP)
\begin{eqnarray}\label{(117)a}
b' = g(x,b), \qquad b(x_0)=b_0.
\end{eqnarray}
Since $b$ is a smooth function, we have that $g(x,b)$ and $\partial_b g(x,b)$ are continuous in some open rectangle 
$$
R=\left\lbrace (x,b): x_1< x<x_2,\hspace*{0.1 cm} y_1< b< y_2\right\rbrace
$$ 
that contains the point $(x_0,b_0)$. Then, by the fundamental existence and uniqueness theorem for ordinary differential equation the IVP \eqref{(117)a} has a unique solution in some closed interval $I = [b_0-\epsilon, b_0+\epsilon]$, where $\epsilon$ is a positive number. Moreover, $x_1$ and $x_2$ has to be chosen so that the strip $x_1< x<x_2$ intersects the domain of the solution of \eqref{th2.2}. Observe that replacing $\phi$ into \eqref{eq} we obtain \eqref{102}. This concludes $(ii)$.

The converse follows from a straightforward computation.

\cqd


\begin{proposition}\label{prop2}
Consider an equation of type \eqref{T} describing pseudospherical surfaces, under the condition \eqref{1}, belonging to the class of equations given by Theorem $\ref{teo7.3}$. There is no local isometric immersion in $\mathbb{R}^3$ of a pseudospherical surface determined by a solution $u$ of the equation, for which the coefficients of the second fundamental form depend on $x, t, z_0, \ldots, z_l, w_1, \ldots, w_m$, $v_1, \ldots, v_n$, where $1\leq l  <\infty$, $1\leq m < \infty$ and $1\leq n < \infty$ are finite, but otherwise arbitrary.
\end{proposition}

\noindent
\textbf{Proof}. Since $\eta\neq 0$ we have $f_{21}\neq 0$, on a open set. If $c+(f_{11}/f_{21})^2a+2f_{11}b/f_{21}=0$ then from Lemma \ref{lemma3.3} the equations \eqref{gauss}, \eqref{3.8} and \eqref{3.9} form an inconsistent system. Therefore, $c+(f_{11}/f_{21})^2a+2f_{11}b/f_{21}\neq 0$ and, from Lemma \ref{lemma3.2}, the coefficients of the second fundamental form of such local isometric immersion are universal, and hence \eqref{3.8} and \eqref{3.9} become
\begin{eqnarray}
f_{11} a_t+f_{21} b_t -f_{12} a_x-f_{22} b_x-2b\Delta_{13}+(a-c)\Delta_{23} = 0,&&\label{a7.11}\\
f_{11} b_t+f_{21} c_t -f_{12} b_x-f_{22} c_x+(a-c)\Delta_{13}+2b\Delta_{23} = 0,&&\label{a7.12}
\end{eqnarray}
where $\Delta_{13} =\lambda(m_1\mu - \eta)z_1 $ and $\Delta_{23}=-\lambda m_1z_1$. Hence, since $f_{ij}$ are given by \eqref{th2.3a}, it follows from \eqref{a7.11} and \eqref{a7.12} that, respectively,
\begin{eqnarray}
[a_t + \mu b_t + \lambda(a_x+\mu b_x)z_0]h + \lambda [ m_2(a_x+\mu b_x)-2(m_1\mu -\eta)b - m_1(a-c)]z_1+\eta(b_t+\lambda z_0b_x) = 0,\label{a7.13}
\end{eqnarray}
\begin{eqnarray}
[b_{t} + \mu c_{t} + \lambda (b_{x}+\mu c_{x})z_{0}]h + \lambda [ m_2(b_x+\mu c_x)+(m_1\mu -\eta)(a-c) - 2m_1b]z_1+\eta(c_t+\lambda z_0c_x) = 0.\label{a7.14}
\end{eqnarray}
Differentiating \eqref{a7.13} and \eqref{a7.14} with respect to $z_2$, since $h'\neq 0$, we have
\begin{eqnarray}
a_t + \mu b_t + \lambda(a_x+\mu b_x)z_0 = 0,&&\label{a7.15}\\
b_{t} + \mu c_{t} + \lambda (b_{x}+\mu c_{x})z_{0}=0.&&\label{a7.16}
\end{eqnarray}
Differentiating \eqref{a7.15} and \eqref{a7.16} with respect to $z_0$, since $\lambda\neq 0$, and replacing the result back into \eqref{a7.15} and \eqref{a7.16} leads to
\begin{eqnarray}
a_t + \mu b_t=0, \quad   a_x+\mu b_x = 0,&&\label{a7.17}\\
b_{t} + \mu c_{t}=0, \quad  b_{x}+\mu c_{x}=0.&&\label{a7.18}
\end{eqnarray}
Substituting \eqref{a7.17} and \eqref{a7.18} into \eqref{a7.13} and \eqref{a7.14} and taking the $z_1$ derivative of the remaining expression we get
\begin{eqnarray}
 -2(m_1\mu -\eta)b - m_1(a-c) = 0,&&\label{a7.13.1}\\
 (m_1\mu -\eta)(a-c) - 2m_1b = 0.&&\label{a7.14.1}
\end{eqnarray}
Since the Gauss equation \eqref{gauss} needs to be satisfied we have $(a-c)^2+b^2\neq 0$. Hence, from \eqref{a7.13.1} and \eqref{a7.14.1} we obtain $(m_1\mu-\eta)^2+m_1^2 = 0$, i.e., $m_1=\eta = 0$, which gives a contradiction since $\eta \neq 0$.

\cqd


\begin{proposition}\label{prop3}
Consider an equation of type \eqref{T} describing pseudospherical surfaces, under the condition \eqref{1}, given by the Theorem $\ref{teo7.4}$. There exists a local isometric immersion in $\mathbb{R}^3$ of a pseudospherical surface, defined by a solution $u$, for which the coefficients $a$, $b$ and $c$ of the second fundamental form depend on $x, t, z_0, \ldots, z_l, w_1, \ldots, w_m$, $v_1, \ldots, v_n$, where $1\leq l  <\infty$, $1\leq m < \infty$ and $1\leq n < \infty$ are finite, but otherwise arbitrary, if, and only if,  

\vspace*{0.4 cm}
\noindent
$(i)$ $\mu = m_1 = 0$, $m_2\neq 0$ and $a$, $b$ and $c$ depend only on $t$ and are given by
\begin{eqnarray}\label{d7.1}
a = \pm \sqrt{L(t)}, \quad b=\beta e^{\pm 2m_2t}, \quad c = a \mp \frac{a_t}{m_2},
\end{eqnarray}
where $L(t) = \sigma e^{\pm 2m_2t}-\beta^2 e^{\pm 4m_2t}-1$, with $\sigma$, $\beta$ $\in$ $\mathbb{R}$, $\sigma >0$ and $\sigma^2 >4\beta^2$. The coefficients $a$, $b$, $c$ are defined on a trip of $\mathbb{R}$ where
\begin{eqnarray}\label{d7.2}
log\sqrt{\frac{\sigma -\sqrt{\sigma^2 -4\beta^2}}{2\beta^2}}< \pm m_2t <log\sqrt{\frac{\sigma +\sqrt{\sigma^2 -4\beta^2}}{2\beta^2}}.
\end{eqnarray}
Moreover, the constants $\beta$ and $\sigma$ have to be chosen so that the strip intersects the domain of the solution of \eqref{th2.4}.

\vspace*{0.3 cm}
\noindent
or
\vspace*{0.3 cm}

\noindent
$(ii)$ $\mu=0$, $m_1\neq 0$, $\lambda^2+m_2^2 \neq 0$ and $a$, $b$ and $c$ are functions of $m_1x+m_2t$ and given by
\begin{eqnarray}\label{ab5.1}
a = \pm \sqrt{L(m_1 x+m_2 t)},\quad b = -\beta e^{\pm2 (m_1x+m_2t)} ,\quad c = a\mp a',
\end{eqnarray}
where $L(m_1 x+m_2 t) = \sigma e^{\pm 2(m_1 x+m_2 t)}-\beta^2 e^{\pm 4(m_1 x+m_2 t)}-1$, with $\sigma$, $\beta$ $\in$ $\mathbb{R}$, $\sigma >0$ and $\sigma^2 >4\beta^2$. The coefficients $a$, $b$, $c$ are defined on a trip of $\mathbb{R}$ where
\begin{eqnarray}\label{d7.2.1}
log\sqrt{\frac{\sigma -\sqrt{\sigma^2 -4\beta^2}}{2\beta^2}}< \pm (m_1 x+ m_2 t) <log\sqrt{\frac{\sigma +\sqrt{\sigma^2 -4\beta^2}}{2\beta^2}}.
\end{eqnarray}
Moreover, the constants $\beta$ and $\sigma$ have to be chosen so that the strip intersects the domain of the solution of \eqref{th2.4}.

\vspace*{0.3 cm}
\noindent
or
\vspace*{0.3 cm}

\noindent
$(iii)$ $\mu\neq 0$, $(\lambda m_1)^2+m_2^2 \neq 0$ and $a$, $b$ and $c$ are differentiable functions of $m_1x+m_2t$ and given by
\begin{eqnarray}\label{129.a}
 a &=& \frac{1}{2\mu}[\pm \mu\sqrt{\Delta}-(\mu^2-1)b+\beta e^{\pm 2(m_1 x+m_2 t)}],\nonumber\\
 c &=& \frac{1}{2\mu}[\pm \mu\sqrt{\Delta}+(\mu^2-1)b-\beta e^{\pm 2(m_1 x+m_2 t)}],\\
 \Delta &=& \frac{[(\mu^2-1)b-\beta e^{\pm 2(m_1 x+m_2 t)}]^2 -4\mu^2(1-b^2)}{\mu^2}> 0\nonumber
\end{eqnarray}
where $b$ satisfies the ordinary differential equation
\begin{eqnarray}\label{102.1}
[\mu(1+\mu^2)\sqrt{\Delta}\pm (\mu^2+1)^2b\mp (\mu^2-1)\beta e^{\pm 2(m_1 x+m_2 t)}]b'\hspace*{3cm}\nonumber\\
+ 2 [\mp\mu (1+\mu^2)\sqrt{\Delta}- \beta(\mu^2-1)e^{\pm 2(m_1 x+m_2 t)}]b+ 2\beta^2 e^{\pm 4(m_1 x+m_2 t)}=0
\end{eqnarray}
\end{proposition}

\noindent
\textbf{Proof}. If $f_{21} \equiv 0$ then $\mu=m_1 =0$ and $m_2\neq 0$. From Lemma \ref{lemma3.4} the coefficients of the second fundamental form of such local isometric immersion are universal, and hence \eqref{3.8} and \eqref{3.9} become
\begin{eqnarray}
[a_t +\lambda a_x z_0\mp m_2(a -c)]h-a_x\psi - m_2b_x = 0,\quad [b_t+\lambda b_x z_0\mp  2m_2 b]h-b_x\psi - m_2c_x = 0.\label{d5.4}
\end{eqnarray}
Differentiating \eqref{d5.4} with respect to $z_2$, since $h'\neq 0$ on a open set, we obtain
\begin{eqnarray}
a_t +\lambda a_x z_0\mp m_2(a -c)=0, \quad b_t+\lambda b_x z_0\mp  2m_2 b=0.\label{d5.5}
\end{eqnarray}
Differentiating \eqref{d5.5} with respect to $z_0$ and replacing the result back into \eqref{d5.5}, we get
\begin{eqnarray}
\lambda a_x = \lambda b_x = 0,\quad a_t \mp m_2(a-c) = 0, \quad b_t \mp 2m_2b = 0.\label{d5.6}
\end{eqnarray}
Substituting \eqref{d5.6} into \eqref{d5.4} we finally have
\begin{eqnarray}\label{d5.7}
a_x \psi + m_2b_x = 0, \quad b_x\psi +m_2c_x = 0.
\end{eqnarray}
Taking the derivative of the Gauss equation \eqref{gauss} with respect to $x$ leads to $a_xc+ac_x - 2bb_x = 0$. Replacing \eqref{d5.7} in the latter, we have
\begin{eqnarray}\label{(136)}
a_x \left[c + \left(\frac{\psi}{m_2}\right)^2a + 2\frac{\psi}{m_2} b\right] = 0.
\end{eqnarray}
If $a_x\neq 0$ then differentiating the first equation in \eqref{d5.7} with respect to $z_0$ and $z_1$ gives us $\psi_{,z_0} = \psi_{,z_1}=0$ and, thus, $\psi = \alpha m_2$, where $\alpha$ denotes a arbitrary constant. From \eqref{(136)} we can see that $\psi\neq 0$, since $c\neq 0$, and
\begin{eqnarray}\label{d5.8}
b = \pm 1 -\alpha a, \quad c = \alpha^2 a\mp 2\alpha.
\end{eqnarray}
Substituting \eqref{d5.8} into \eqref{d5.6} leads to
$$
a_t \mp m_2(a-\alpha^2 a\pm 2\alpha) = 0, \quad -\alpha a_t \mp 2m_2(\pm 1-\alpha a) = 0.
$$
In the above equations, adding the second to the first multiplied by $\alpha$ leads to $a=\pm 2/\alpha$, which replaced in the first equation gives us $m_2 = 0$ and, thus, a contradiction since $m_2\neq 0$.

Therefore, $a_x = 0$ and by \eqref{d5.7} we have $b_x = c_x = 0$. Thus, $a$, $b$ and $c$ depend only on $t$. It follows from \eqref{d5.6} that
\begin{eqnarray}\label{d5.9}
b = \beta e^{\pm 2m_2t}, \quad c = a\mp \frac{a_t}{m_2},
\end{eqnarray}
where $\beta$ is a constant. Replacing \eqref{d5.9} into the Gauss equation leads to $a= \pm \sqrt{L(t)}$ where $L(t) = \sigma e^{\pm 2m_2t} - \beta^2 e^{\pm 4m_2t} - 1$, $\sigma>0$ is a constant and $\sigma^2>4\beta^2$. This $a$ together with \eqref{d5.9} gives \eqref{d7.1}, where $a$ is defined on the trip described by \eqref{d7.2}. Observe that $\psi$ and $\lambda$ are still arbitrary. This concludes $(i)$.

Suppose $f_{21}\neq 0$ on a non-empty open set. If $c+(f_{11}/f_{21})^2a+2f_{11}b/f_{21}=0$ then from Lemma \ref{lemma3.3} the equations \eqref{gauss}, \eqref{3.8} and \eqref{3.9} form an inconsistent system. Therefore, we have $c+(f_{11}/f_{21})^2a+2f_{11}b/f_{21}\neq 0$ and, from Lemma \ref{lemma3.2}, the coefficients of the second fundamental form of such local isometric immersion are universal, and hence \eqref{3.8} and \eqref{3.9} become
\begin{eqnarray}
f_{11} a_t+f_{21} b_t -f_{12} a_x-f_{22} b_x-2b\Delta_{13}+(a-c)\Delta_{23} = 0,&&\label{7.21}\\
f_{11} b_t+f_{21} c_t -f_{12} b_x-f_{22} c_x+(a-c)\Delta_{13}+2b\Delta_{23} = 0,&&\label{7.22}
\end{eqnarray}
where $\Delta_{13} = \pm \mu(\lambda m_1z_0+m_2)h\mp \mu m_1\psi$ and $\Delta_{23}=\mp (m_2+\lambda m_1z_0)h \pm m_1\psi$. Hence, since $f_{ij}$ are given by \eqref{th2.4a}, it follows from \eqref{7.21} and \eqref{7.22} that, respectively,
\begin{eqnarray}
[a_t+\mu b_t +\lambda (a_x+\mu b_x)z_0\mp (a +2\mu b-c)(\lambda m_1z_0 +m_2)]h-[a_x+\mu b_x \mp m_1(a+2\mu b-c)]\psi\nonumber\\
+\sqrt{1+\mu^2}(m_1b_t - m_2b_x) = 0,\label{7.23}
\end{eqnarray}
\begin{eqnarray}
[b_t+\mu c_t +\lambda (b_x+\mu c_x)z_0\pm (\mu a -2 b-\mu c)(\lambda m_1z_0 +m_2)]h-[b_x+\mu c_x \pm m_1(\mu a-2 b-\mu c)]\psi\nonumber\\
+\sqrt{1+\mu^2}(m_1c_t - m_2c_x) = 0.\label{7.24}
\end{eqnarray}
Differentiating \eqref{7.23} and \eqref{7.24} with respect to $z_2$ leads, since $h'\neq 0$ on a open set, to
\begin{eqnarray}
a_t+\mu b_t +\lambda (a_x+\mu b_x)z_0\mp (a +2\mu b-c)(\lambda m_1z_0 +m_2)=0,&&\label{7.25}\\
b_t+\mu c_t +\lambda (b_x+\mu c_x)z_0\pm (\mu a -2 b-\mu c)(\lambda m_1z_0 +m_2)=0.&&\label{7.26}
\end{eqnarray}
Differentiating \eqref{7.25} and \eqref{7.26} with respect to $z_0$ and replacing the result back into the latter two equations we get
\begin{eqnarray}\label{7.27}
a_t+\mu b_t\mp m_2(a +2\mu b-c)=0, \quad b_t+\mu c_t \pm m_2(\mu a -2 b-\mu c)=0, 
\end{eqnarray}
and
\begin{eqnarray}\label{7.28}
\lambda[a_x+\mu b_x\mp m_1(a +2\mu b-c)]=0, \quad \lambda[b_x+\mu c_x\pm m_1(\mu a -2 b-\mu c)]=0
\end{eqnarray}

Finally, substituting \eqref{7.27} and \eqref{7.28} back into \eqref{7.23} and \eqref{7.24}, we obtain
\begin{eqnarray}
-[a_x+\mu b_x \mp m_1(a+2\mu b-c)]\psi +\sqrt{1+\mu^2}(m_1b_t - m_2b_x) = 0,&&\label{143}\\
-[b_x+\mu c_x \pm m_1(\mu a-2 b-\mu c)]\psi +\sqrt{1+\mu^2}(m_1c_t - m_2c_x) = 0.&&\label{144}
\end{eqnarray}
Multiplying the first equation in \eqref{7.27} by $m_2$ and adding the result to the first equation in \eqref{7.28} multiplies by $\lambda m_1$ leads to
\begin{eqnarray}\label{145}
a+2\mu b - c = \pm \frac{1}{M}[m_2(a_t+\mu b_t)+\lambda^2 m_1(a_x+\mu b_x)],
\end{eqnarray}
and the same operation with the second equation of \eqref{7.27} and \eqref{7.28} leads to
\begin{eqnarray}\label{146}
\mu a -2 b -\mu c = \mp \frac{1}{M}[m_2(b_t+\mu c_t)+\lambda^2 m_1(b_x+\mu c_x)],
\end{eqnarray}
where $M = (\lambda m_1)^2+m_2^2$ is a nonzero constant. Replacing \eqref{145} and \eqref{146} into \eqref{143} and \eqref{144}, we obtain
\begin{eqnarray}
m_2 \psi (m_1a_t - m_2a_x) + (\mu m_2\psi +M\sqrt{1+\mu^2})(m_1b_t-m_2b_x) = 0,&&\label{147}\\
m_2 \psi (m_1b_t - m_2b_x) + (\mu m_2\psi +M\sqrt{1+\mu^2})(m_1c_t-m_2c_x) = 0.&&\label{148}
\end{eqnarray}
Differentiating the Gauss equation \eqref{gauss} with respect to $t$ and multiplying the result by $m_1$ and doing the same thing with $x$ and $m_2$, we get
\begin{eqnarray*}
m_1 a_t c+m_1ac_t - 2m_1bb_t = 0,&&\\
m_2 a_x c+m_2ac_x - 2m_2bb_x = 0.&&
\end{eqnarray*}
From the two latter equations we obtain
\begin{eqnarray}\label{G}
(m_1a_t -m_2 a_x)c +  (m_1c_t -m_2 c_x)a - 2b(m_1b_t -m_2 b_x)=0.
\end{eqnarray}

Suppose $m_2\psi \neq 0$. Replacing \eqref{147} and \eqref{148} in \eqref{G} leads to
\begin{eqnarray}\label{149}
(m_1c_t -m_2 c_x)[a+Q^2c + 2Qb] = 0, \quad Q = \frac{\mu m_2\psi + M\sqrt{1+\mu^2}}{m_2\psi}.
\end{eqnarray}
If $m_1c_t -m_2 c_x\neq 0$, then $\psi$ is a constant and, by \eqref{149} and the Gauss equation, we have
\begin{eqnarray}\label{150}
a = Q^2 c \mp 2Q, \quad b = \pm 1 -Qc.
\end{eqnarray}
Since $a\neq 0$ we have $Q\neq 0$. Substituting \eqref{150} into \eqref{7.27}, we get
$$
(Q-\mu )Qc_t \mp m_2[Q^2 c\mp 2Q - c +2\mu (\pm 1-Qc)] = 0, \quad -(Q-\mu )c_t \pm m_2[\mu(Q^2 c\mp 2Q - c) -2 (\pm 1-Qc)] = 0.
$$

In the latter equations, adding the second multiplied by $Q$ to the first gives us $c=constant$, which implies from \eqref{150} that $a$ and $b$ are constants. But, $a$, $b$ and $c$ constants imply from \eqref{7.27} that $a-c=0$ and $b=0$, which contradicts the Gauss equation \eqref{gauss}. Therefore $m_1c_t -m_2 c_x = 0$ and, thus, from \eqref{147} and \eqref{148} we have $m_1b_t-m_2b_x=0$ and $m_1a_t-m_2a_x = 0$.

On the other hand, if $m_2\psi = 0$ then from \eqref{147} and \eqref{148} we have $m_1b_t-m_2b_x=0$ and $m_1c_t -m_2 c_x = 0$, which replaced into \eqref{G} since $c\neq 0$ leads to $m_1a_t-m_2a_x = 0$.

Therefore, for arbitrary $m_2\psi$ we have shown that
\begin{eqnarray}
\begin{array}{ll}\label{7.31}
a = \phi_1(m_1x+m_2t),\quad b = \phi_2(m_1x+m_2t),\quad c = \phi_3(m_1x+m_2t),
\end{array}
\end{eqnarray}
where $\phi_i$, $i=1,2,3$, are real and differentiable functions and, by Lemma \ref{lemma3.1}, $\phi_1\phi_3\neq 0$ on a open set. Replacing \eqref{7.31} into \eqref{7.27} and \eqref{7.28} and observing that $(\lambda m_1)^2+m_2^2\neq 0$, we obtain
\begin{eqnarray}
\phi_1' + \mu \phi_2'\mp (\phi_1+2\mu \phi_2 - \phi_3) = 0,&&\label{b5.12}\\
\phi_2' + \mu \phi_3' \pm (\mu \phi_1 -2\phi_2 - \mu \phi_3) = 0.&&\label{b5.13}
\end{eqnarray}
From \eqref{b5.12} we obtain $\phi_3$ in terms of $\phi_1$, $\phi_1'$, $\phi_2$ and $\phi_2'$, which replaced into \eqref{b5.13} implies that
\begin{eqnarray}\label{b5.14}
\mu\phi_1' = \pm (1+\mu^2)\phi_2 - \mu^2\phi_2' \pm \beta e^{\pm 2(m_1 x+m_2t)}.
\end{eqnarray}
If $\mu=0$, then from \eqref{b5.14} and \eqref{b5.12} we have $b = -\beta e^{\pm 2(m_1 x+m_2t)}$ and $c = a\mp a'$. Using the latter and Gauss equation leads to \eqref{ab5.1}, where $a$ is defined on the trip described by \eqref{d7.2.1}. Observe that $\lambda$ and $\psi$ are still arbitrary. This concludes $(ii)$.

If $\mu \neq 0$, then from \eqref{b5.14} we have $\phi_1'$, which replaced into \eqref{b5.12} implies that
\begin{eqnarray}\label{112.b}
\phi_3 = \phi_1 + \phi(m_1 x+m_2 t), \quad \phi(m_1 x+m_2 t) = \frac{\mu^2-1}{\mu}\phi_2 - \frac{\beta}{\mu}e^{\pm 2(m_1x+m_2t)}.
\end{eqnarray}
Substituting the latter into the Gauss equation we obtain $\phi_1^2+\phi_1\phi(m_1 x+m_2 t)-\phi_2^2 =-1$, which resolved as a second degree equation in terms of $\phi_1$ leads to
$$
\phi_1 = \frac{-\phi(m_1 x+m_2 t)\pm \sqrt{\Delta}}{2}, \quad \Delta = \phi(m_1 x+m_2 t)^2-4[1-\phi_2(m_1 x+m_2 t)^2]> 0.
$$

Hence, using \eqref{112.b} we also have $\phi_3$ in terms of $\phi_2=\phi_2(m_1 x+m_2 t)$ as in \eqref{129.a}, which replaced into \eqref{b5.14} gives us
\begin{eqnarray}\label{eqt}
[(1+\mu^2)\sqrt{\Delta}\pm (\mu^2-1)\phi \pm 4\mu b]b' \mp 2(1+\mu^2)b\sqrt{\Delta}-2\beta e^{\pm 2(m_1 x+m_2 t)}\phi =0.
\end{eqnarray}
Observe that, if the coefficient of $b'$ in \eqref{eqt} vanishes, we have
\begin{eqnarray*}
(1+\mu^2)\sqrt{\Delta}\pm (\mu^2-1)\phi \pm 4\mu b=0,&&\\
\mp 2(1+\mu^2)b\sqrt{\Delta}-2\beta e^{\pm 2(m_1 x+m_2 t)}\phi =0.&&
\end{eqnarray*}
In the latter two equations, replacing $(1+\mu^2)\sqrt{\Delta}$ of the first into the second implies that
\begin{eqnarray*}
0 &=& \mp 2 b[\mp (\mu^2-1)\phi \mp 4\mu b]-2\beta e^{\pm 2(m_1 x+m_2 t)}\phi\\
  &=& 2\phi [(\mu^2-1)b-\beta e^{\pm 2(m_1 x+m_2 t)}]+8\mu b^2\\
  &=&2\phi \mu\phi + 8\mu b^2\\
  &=& 2\mu (\phi^2 +4b^2),
\end{eqnarray*}
and then $\phi^2 +4b^2=0$, since $\mu\neq 0$. However, $\phi^2 +4b^2=0$ if, and only if, $\phi=b=0$ which implies by \eqref{7.12b} that $a=c$. But, $a=c$ and $b=0$ contradict the Gauss equation \eqref{gauss}. Therefore, the coefficient of $b'$ in the equation \eqref{eqt} does not vanish in a non-empty open set. That means we can write $b' = g(x,b)$, where $g$ is a differentiable function defined, from \eqref{eqt}, by
$$
g(x,b) = \frac{\pm 2\eta(1+\mu^2)b\sqrt{\Delta}+2\eta\beta e^{\pm 2\eta x}\phi}{(1+\mu^2)\sqrt{\Delta}\pm (\mu^2-1)\phi \pm 4\mu b}.
$$

Let $x_0$ be an arbitrarily fixed point and consider the following Initial Value Problem (IVP)
\begin{eqnarray}\label{(117)at}
b' = g(x,b), \qquad b(x_0)=b_0.
\end{eqnarray}
Since $b$ is a smooth function, we have that $g(x,b)$ and $\partial_b g(x,b)$ are continuous in some open rectangle 
$$
R=\left\lbrace (x,b): x_1< x<x_2,\hspace*{0.1 cm} y_1< b< y_2\right\rbrace
$$ 
that contains the point $(x_0,b_0)$. Then, by the fundamental existence and uniqueness theorem for ordinary differential equation the IVP \eqref{(117)at} has a unique solution in some closed interval $I = [b_0-\epsilon, b_0+\epsilon]$, where $\epsilon$ is a positive number. Moreover, $x_1$ and $x_2$ has to be chosen so that the strip $x_1< x<x_2$ intersects the domain of the solution of \eqref{th2.4}. Observe that replacing $\phi$ into \eqref{eqt} we obtain \eqref{102.1}. This concludes $(iii)$.

The converse follows from a straightforward computation.

\cqd


\begin{proposition}\label{prop4}
Consider an equation of type \eqref{T} describing pseudospherical surfaces, under the condition \eqref{1}, given by Theorem $\ref{teo7.5}$-$(i)$. There is no local isometric immersion in $\mathbb{R}^3$ of a pseudospherical surface determined by a solution $u$ of the equation, for which the coefficients of the second fundamental form depend on $x, t, z_0, \ldots, z_l, w_1, \ldots, w_m$, $v_1, \ldots, v_n$, where $1\leq l  <\infty$, $1\leq m < \infty$ and $1\leq n < \infty$ are finite, but otherwise arbitrary.
\end{proposition}

\noindent
\textbf{Proof}. Since $\eta\neq 0$ we can not have $f_{21}=0$ on a open set. Therefore $f_{21}\neq 0$ on a open set. 

If $c+(f_{11}/f_{21})^2a+2f_{11}b/f_{21}=0$ then from Lemma \ref{lemma3.3} the equations \eqref{gauss}, \eqref{3.8} and \eqref{3.9} form an inconsistent system. Hence, we have $c+(f_{11}/f_{21})^2a+2f_{11}b/f_{21}\neq 0$ and, from Lemma \ref{lemma3.2}, the coefficients of the second fundamental form of such local isometric immersion are universal, and hence \eqref{3.8} and \eqref{3.9} become
\begin{eqnarray}
f_{11} a_t+f_{21} b_t -f_{12} a_x-f_{22} b_x-2b\Delta_{13}+(a-c)\Delta_{23} = 0,&&\label{7.21a}\\
f_{11} b_t+f_{21} c_t -f_{12} b_x-f_{22} c_x+(a-c)\Delta_{13}+2b\Delta_{23} = 0,&&\label{7.22a}
\end{eqnarray}
where 
\begin{eqnarray}\label{7.34}
\begin{array}{ll}
\Delta_{13} = \left[\mu \left(m-\frac{b\tau}{a}\right)-\frac{\tau\eta}{a}\right](\mp \phi_{12}+\tau \varphi e^{\pm \tau z_1}f_{11}),\quad 
\Delta_{23}= -\left(m-\frac{b\tau}{a}\right)(\mp \phi_{12}+\tau \varphi e^{\pm \tau z_1}f_{11}),
\end{array}
\end{eqnarray}
with $\phi_{12} = [\pm \tau(az_0+b)\varphi +az_1\varphi']e^{\pm \tau z_1}\mp \lambda az_1/\tau$. It follows from \eqref{7.21a} and \eqref{7.22a} that, respectively,
\begin{eqnarray}\label{7.35}
[a_t+\mu b_t+\lambda (a_x+\mu b_x)z_0]f_{11}-(a_x+\mu b_x)\phi_{12}+\eta b_t + \eta(\lambda z_0\mp \tau e^{\pm \tau z_1}\varphi)b_x-2b\Delta_{13}+(a-c)\Delta_{23} = 0,
\end{eqnarray}
\begin{eqnarray}\label{7.36}
[b_t+\mu c_t+\lambda (b_x+\mu c_x)z_0]f_{11}-(b_x+\mu c_x)\phi_{12}+\eta c_t + \eta(\lambda z_0\mp \tau e^{\pm \tau z_1}\varphi)c_x+(a-c)\Delta_{13}+2b\Delta_{23} = 0.
\end{eqnarray}
Differentiating \eqref{7.35} and \eqref{7.36} with respect to $z_2$ and using \eqref{7.34}, since $f_{11,z_2}\neq 0$, we get
\begin{eqnarray}\label{7.37}
a_t+\mu b_t+\lambda (a_x+\mu b_x)z_0-\left\lbrace 2b\left[\mu \left(m-\frac{b\tau}{a}\right)-\frac{\tau\eta}{a}\right]+(a-c)\left(m-\frac{b\tau}{a}\right) \right\rbrace \tau \varphi e^{\pm \tau z_1}  = 0,
\end{eqnarray}
\begin{eqnarray}\label{7.38}
b_t+\mu c_t+\lambda (b_x+\mu c_x)z_0+\left\lbrace (a-c)\left[\mu \left(m-\frac{b\tau}{a}\right)-\frac{\tau\eta}{a}\right]-2b\left(m-\frac{b\tau}{a}\right) \right\rbrace \tau \varphi e^{\pm \tau z_1} = 0.
\end{eqnarray}
Differentiating \eqref{7.37} and \eqref{7.38} with respect to $z_1$ and observing that $\tau\varphi \neq 0$, we have
\begin{eqnarray}\label{7.39}
\begin{array}{rr}
2b\left[\mu \left(m-\frac{b\tau}{a}\right)-\frac{\tau\eta}{a}\right]+(a-c)\left(m-\frac{b\tau}{a}\right) = 0,\\
-2b\left(m-\frac{b\tau}{a}\right)+(a-c)\left[\mu \left(m-\frac{b\tau}{a}\right)-\frac{\tau\eta}{a}\right]= 0.
\end{array}
\end{eqnarray}
Since $b^2+(a-c)^2\neq 0$, it follows from \eqref{7.39} that
$$
\left[\mu \left(m-\frac{b\tau}{a}\right)-\frac{\tau\eta}{a}\right]^2+\left(m-\frac{b\tau}{a}\right)^2 = 0,
$$
which implies $\Delta_{13} = \Delta_{23} = 0$ and, thus, contradicts \eqref{3.10a}.

\cqd


\begin{proposition}\label{prop5}
Consider an equation of type \eqref{T} describing pseudospherical surfaces, under the condition \eqref{1}, given by Theorem $\ref{teo7.5}$-$(ii)$. There is no local isometric immersion in $\mathbb{R}^3$ of a pseudospherical surface determined by a solution $u$ of the equation, for which the coefficients of the second fundamental form depend on $x, t, z_0, \ldots, z_l, w_1, \ldots, w_m$, $v_1, \ldots, v_n$, where $1\leq l  <\infty$, $1\leq m < \infty$ and $1\leq n < \infty$ are finite, but otherwise arbitrary.
\end{proposition}

\noindent
\textbf{Proof}. Suppose $f_{21}\neq 0$ on a open set. If $c+(f_{11}/f_{21})^2a+2f_{11}b/f_{21}=0$ then from Lemma \ref{lemma3.3} the equations \eqref{gauss}, \eqref{3.8} and \eqref{3.9} form an inconsistent system. Hence, we have $c+(f_{11}/f_{21})^2a+2f_{11}b/f_{21}\neq 0$ and, from Lemma \ref{lemma3.2}, the coefficients of the second fundamental form of such local isometric immersion are universal, and hence \eqref{3.8} and \eqref{3.9} become
\begin{eqnarray}
f_{11} a_t+f_{21} b_t -f_{12} a_x-f_{22} b_x-2b\Delta_{13}+(a-c)\Delta_{23} = 0,&&\label{11.1}\\
f_{11} b_t+f_{21} c_t -f_{12} b_x-f_{22} c_x+(a-c)\Delta_{13}+2b\Delta_{23} = 0,&&\label{11.2}
\end{eqnarray}
where by \eqref{(7.2)} we have
\begin{eqnarray}\label{11.3}
\begin{array}{ll}
\Delta_{13} = (\phi_{32}\mp \sqrt{1+\mu^2}\phi_{12})f_{11}\mp \frac{\theta +a\mu\eta}{a\sqrt{1+\mu^2}}\phi_{12},\\
\Delta_{23} = (\mu \phi_{32}\mp \sqrt{1+\mu^2}\phi_{22})f_{11}+\eta \phi_{32}\mp \frac{\theta +a\mu\eta}{a\sqrt{1+\mu^2}}\phi_{22}.
\end{array}
\end{eqnarray}
Differentiating \eqref{11.1} and \eqref{11.2} with respect to $z_2$, we obtain since $f_{11,z_2}\neq 0$, respectively,
\begin{eqnarray}
a_t+\mu b_t+\lambda (a_x+\mu b_x)z_0 - 2b(\phi_{32}\mp\sqrt{1+\mu^2}\phi_{12})+(a-c)(\mu\phi_{32}\mp \sqrt{1+\mu^2}\phi_{22}) = 0,\label{a11.4}\\
b_t+\mu c_t+\lambda (b_x+\mu c_x)z_0 +(a-c)(\phi_{32}\mp\sqrt{1+\mu^2}\phi_{12})+2b(\mu\phi_{32}\mp \sqrt{1+\mu^2}\phi_{22}) = 0.\label{a11.5}
\end{eqnarray}
Differentiating \eqref{a11.4} and \eqref{a11.5} with respect to $z_1$, since $b^2+(a-c)^2\neq 0$, we conclude that
$$
(\phi_{32}\mp\sqrt{1+\mu^2}\phi_{12})_{,z_1} = (\mu\phi_{32}\mp \sqrt{1+\mu^2}\phi_{22})_{,z_1}=0
$$
if, and only if, $m_1\theta e^{\theta z_0}-\lambda = 0$, i.e., if and only if $\lambda = m_1=0$, which implies $\Delta_{12} = 0$ and, thus, a contradiction.

On the other hand, if $f_{21} = 0$, on a open set, then we have $\mu=\eta = 0$. It follows from Lemma \ref{lemma3.4} that the coefficients of the second fundamental form of such local isometric immersion are universal, and hence \eqref{3.8} and \eqref{3.9} become
\begin{eqnarray}
f_{11} a_t -f_{12} a_x-f_{22} b_x-2b\Delta_{13}+(a-c)\Delta_{23} = 0,&&\label{11.6}\\
f_{11} b_t -f_{12} b_x-f_{22} c_x+(a-c)\Delta_{13}+2b\Delta_{23} = 0,&&\label{11.7}
\end{eqnarray}
where by \eqref{(7.2)} we have $\Delta_{13} = (\phi_{32}\mp \phi_{12})f_{11}\mp \theta \phi_{12}/a$ and $\Delta_{23} = \mp\phi_{22}f_{11}\mp \theta\phi_{22}/a$. Differentiating \eqref{11.6} and \eqref{11.7} with respect to $z_2$, since $f_{11,z_2}\neq 0$, we have, respectively,
\begin{eqnarray}
a_t+\lambda a_xz_0 - 2b(\phi_{32}\mp\phi_{12})\mp(a-c) \phi_{22} = 0,\label{11.4}\\
b_t+\lambda b_xz_0 +(a-c)(\phi_{32}\mp\phi_{12})\mp 2b \phi_{22} = 0,\label{11.5}
\end{eqnarray}
where $\phi_{32}\mp\phi_{12} = \pm (m_1\theta e^{\theta z_0}-\lambda)/a$ and $\phi_{22}=\mp (m_1\theta e^{\theta z_0}-\lambda)z_1$. Differentiating \eqref{11.4} and \eqref{11.5} with respect to $z_1$, since $m_1\theta e^{\theta z_0}-\lambda\neq 0$, we have $b=a-c=0$, which contradicts the Gauss equation.

\cqd

Finally, the proof of Theorem \ref{teo} follows from Propositions \ref{prop1}-\ref{prop5}.

\cqd

\noindent
\textbf{Acknowledgments:} We are grateful to the Department of Mathematics and Statistics of the McGill University for its hospitality while this paper was being prepared. We express our sincere gratitude to Keti Tenenblat for invaluable suggestions. This work was supported by Ministry of Science and Technology, Brazil, CNPq Proc. 248877/2013-5  and by NSERC Grant RGPIN 105490-2011.


\vspace{1 cm}
\noindent
Tarc\'isio Castro Silva\\
Department of Mathematics and Statistics, McGill University, Canada\\
e-mail: tarcisio.castrosilva@mail.mcgill.ca

\vspace{0.4 cm}
\noindent
Niky Kamran\\
Department of Mathematics and Statistics, McGill University, Canada\\
e-mail: nkamran@math.mcgill.ca

\end{document}